\documentclass{article}

\newcommand{\R}{\mathbb{R}}
\newcommand{\C}{\mathbb{C}}

\newcommand{\at}{{\mathcal A}_{\bf T}}
\newcommand{\iit}{{\mathcal I}_{\bf T}}
\newcommand{\LL}{L^{1}_{\omega_{\bf T}}}
\newcommand{\M}{{\mathcal M}_{\omega_{\bf T}}}

\newcommand{\mcu}{{\mathcal U}}
\newcommand{\qm}{{\mathcal Q}{\mathcal M}}

\newcommand{\jjt}{{\mathcal J}_{\bf T}}
\newcommand{\lt}{{\mathcal L}_{\bf T}}
\newcommand{\dt}{{\mathcal U}_{\bf T}}
\newcommand{\Uu}{\mathcal U}

\newcommand{\sia}{\widehat{\iit}}

\usepackage[ansinew]{inputenc}
\usepackage[T1]{fontenc}
\usepackage[english]{babel}
\usepackage{graphicx}
\usepackage{color}
\usepackage{amsfonts}
 \newtheorem{thm}{Theorem}[section]
 \newtheorem{cor}[thm]{Corollary}
 \newtheorem{lem}[thm]{Lemma}
 \newtheorem{prop}[thm]{Proposition}
  \newtheorem{defn}[thm]{Definition}
  
\begin{document}
%
%
%
%
%
%
%
%
%

\title
 {On the generation of Arveson weakly continuous semigroups}
\author{Jean Esterle}

\maketitle


\begin{abstract}
We consider here one-parameter semigroups ${\bf T}=(T(t))_{t>0}$ of bounded operators on a Banach space $X$ which are weakly continuous in the sense of Arveson. For such a semigroup ${\bf T}$ denote by ${\mathcal M}_{\omega_{\bf T}}$ the convolution algebra consisting in those measures $\mu$ on $(0,+\infty)$ such that $\int_0^{+\infty}\Vert T(t)\Vert d\vert \mu \vert (t)<+\infty.$ The Pettis integral $\int_0^{+\infty}T(t)d\mu(t)$ defines for $\mu \in {\mathcal M}_{\omega_{\bf T}}$ a bounded operator $\phi_{\bf T}(\mu)$ on $X.$ Identifying
the space $L^1_{\omega_{\bf T}}$ of (classes of) measurable functions $f$ satisfying $\int_0^{+\infty}\vert f(t)\Vert T(t)\Vert dt< +\infty$ to a closed subspace ${\mathcal M}_{\omega_{\bf T}}$ in the usual way, we define the Arveson ideal $\iit$
of the semigroup to be the closure in ${\mathcal B}(X)$ of $\phi_{\bf T}(L^1_{\omega_{\bf T}}).$ Using a variant of a procedure introduced a long time ago by the author we introduce a dense ideal $\dt$ of $\iit,$ which is a Banach algebra with respect to a suitable norm $\Vert .\Vert_{\dt},$ such that $\lim \sup_{t\to 0^+}\Vert T(t)\Vert_{{\mathcal B}(\dt)}<+\infty.$ The normalized Arveson ideal $\jjt$ is the closure of $\iit$ in ${\mathcal B}(\dt).$ The Banach algebra $\jjt$ has a sequential approximate identity and is isometrically isomorphic to a closed ideal of its multiplier algebra ${\mathcal M}(\jjt).$ The Banach algebras $\dt,$ $\iit$ and $\jjt$ are "similar", and the map $S_{u/v}\to S_{au/av}$ defines when $a$ generates a dense principal ideal of $\dt$ a pseudobounded isomorphism from the algebre $\qm(\jjt)$ of quasimultipliers on $\jjt$ onto the quasimultipliers algebras $\qm(\dt)$ and $\qm(\iit).$ 

We define the generator $A_{\bf T}$ of the semigroup $\bf T$ to be a quasimultiplier on $\iit,$ or ,equivalently, on $\jjt.$ Every character $\chi$ on $\iit$ has an extension $\tilde \chi$ to $\qm(\iit).$ Let $Res_{ar} (A_{\bf T})$ be the complement of the set $\{ \tilde \chi (A_{\bf T})\}_{\chi \in \sia}.$ The quasimultiplier $A-\mu I$ has an inverse belonging to $\jjt$ for $\mu \in Res_{ar} (A_{\bf T}),$ which allows to consider this inverse as a "regular" quasimultiplier on the Arveson ideal $\iit.$ The usual resolvent formula holds in this context for $Re(\mu)>
\lim_{t\to +\infty}{log \Vert T(t)\Vert\over t}.$

Set $\Pi_{\alpha}^+:=\{ z \in \C \ | \ Re(z) >\alpha\}.$ We revisit the functional calculus associated to the generator $A_{\bf T}$ by defining $F(-A_{\bf T})\in \jjt$ by a Cauchy integral when $F$ belongs to the Hardy space $H^1(\Pi_{\alpha}^+)$ for some $\alpha < -\lim_{t\to +\infty} {log\Vert T(t)\vert\over t}.$ We then define $F(-A_{\bf T})$ as a quasimultiplier on $\jjt$ and $\iit$ when $F$ belongs to the Smirnov class on $\Pi_{\alpha}^+,$ and $F(-A_{\bf T})$ is a regular quasimultiplier on $\jjt$ and $\iit$ if $F$ is bounded on $\Pi_{\alpha}^+.$ If $F(z)=e^{-zt}$ for some $t>0,$ then $F(-A_{\bf T})=T(t),$ and if $F(z)=-z,$ we indeed have $F(-A_{\bf T})=A_{\bf T}.$
\end{abstract}

keywords: semigroup of bounded operators, Arveson pairs, Arveson spectrum, Pettis integral, infinitesimal generator, resolvent, Laplace transform, holomorphic functional calculus



\maketitle

\section{Introduction}
We study here the generator of a semigroup ${\bf T}=(T(t))_{t>0}$ of bounded operators on a Banach space $X$ which is weakly continuous with respect to a dual pair $(X,X_*)$ satisfying the conditions introduced by Arveson in 1974 in his seminal paper \cite{ar} on group representations. Semigroups satisfying some weaker properties have been considered in 1953 by Feller, who studied semigroups ${\bf T}=(T(t))_{t>0}$ of bounded operators on a Banach space $X$ such there exists a nonzero continuous linear form $l$ on $X$ for which the map $t \to <T(t)x,l>$ is measurable for every $x \in X,$ and his work involves the subspace $X_*$ of the dual space $X'$ generated by the family $(l\circ T(t))_{t>0}.$ Feller's Pettis integrals allow to construct operators taking values on some subspace $\tilde X$ of $X_*,$ and he shows that there exists a rather large closed subspace $X_1$  of $\tilde X$ such that the map $t \to T(t)x$ is continuous on $(0,+\infty)$ for $x \in X_1.$

Arveson's conditions are more restrictive, but his 
approach allows  to define a norm-decreasing homomorphism $\phi_{\bf T}$ from the convolution Banach algebra ${\mathcal M}_{\omega_{\bf T}}$ into ${\mathcal B}(X)$,without enlarging the Banach space $X.$ Here ${\mathcal M}_{\omega_{\bf T}}$ denotes the algebra of all Borel measures $\mu$ on $(0,+\infty)$ such that $\int_0^{+\infty}\Vert T(t)\Vert d\vert \mu\vert (t)<+\infty.$ 

Denote by $L^1_{\omega_{\bf T}}$ the convolution algebra of (classes) of measurable functions on the half-line satisfying $\int_{0}^{+\infty}\vert f(t)\vert \Vert T(t)\Vert dt<+\infty,$ identified to a closed ideal of
${\mathcal M}_{\omega_{\bf T}}$, and let $\at$ (resp. $\iit$) be the closure of $\phi_{\bf T}(\M)$ (resp $\phi_{\bf T}(\LL)$) in ${\mathcal B}(X)$ with respect to the operator norm. The ideal  $\iit$ of $\at$ is called the Arveson ideal associated to the semigroup ${\bf T},$ see sections 2 and 3.

In section 4 we introduce the infinitesimal generator $A_{{\bf T}, op}$ as a weakly densely defined operator on the space $X_1:=\left [ \cup_{t>0}T(t)(X)\right ]^-,$ following Feller's approach in \cite{f_0}, and we also introduce the resolvent $R_{op}({\bf T},\lambda)=(A_{{\bf T}, op}-\lambda I)^{-1}$ as a weakly densely defined operator on $X_1.$ If $\lim \sup_{t\to 0^+}\Vert T(t)\Vert<+\infty,$ these partially defined operators are closed operators on $X_1.$ In the general case they can be interpreted as closed operators on $X_0:=X_1/\left [ \cap_{t>0}Ker(T(t)\right ],$ which of course equals $X_1$ if $\cap_{t>0}Ker(T(t)=\{0\}.$

In section 5 we pave the way to a more intrinsic aopproach to the infinitesimal generator of a Arveson weakly continuous semigroup ${\bf T}.$ 
It follows from some results of \cite{e1} that $L^1_{\omega_{\bf T}}$ contains a function $g$ such that $g*L^1_{\omega_{\bf T}}$ is dense in $L^1_{\omega_{\bf T}},$ and so $\iit$ possesses dense principal ideals, and it follows from the weak$^*$ continuity of the semigroup that $u\at\neq \{0\}$ for every $u\in \at \setminus \{0\}.$ This allows us to apply to $\iit$ the theory of quasimultipliers developped a long time ago by the author in \cite{e2}. Recall that if ${\mathcal U}$ is a Banach algebra such that $u{\mathcal U}\neq \{0\}$ for $u\in {\mathcal U}\setminus \{0\}$ and such that the set $\Omega({\mathcal U}):=\{ u \in {\mathcal U} \ | \ [u{\mathcal U}]^-= {\mathcal U}\}$ is nonempty, a quasimultiplier on $\mathcal U$ is a closed densely defined operator on $\mathcal U$ of the form $S=S_{u/v},$ where $u \in \Uu, v \in \Omega(\Uu),$ and where the domain ${\mathcal D}_{S}$ of $S$ is the set of all $x \in \Uu$ such that $ux \in v\Uu.$ If $x \in {\mathcal D},$ then $Sx$ is the unique $y \in \Uu$ such that $ux=vy.$ A subset $\Delta$ of the algebra ${\mathcal {QM}}(\Uu)$ is said to be pseudobounded if $\cap_{S \in \Delta}{\mathcal D}_s$ contains some $u\in \Omega(\Uu)$ such that $\sup_{S\in \Delta}\Vert Su\Vert<+\infty.$ A quasimultiplier $S\in {\mathcal{QM}}(\Uu)$ is said to be regular if there exists $\lambda >0$ such that the set $\{\lambda^nS^n\}_{n\ge 1}$ is pseudobounded, and the family ${\mathcal {QM}}_r(\Uu)$ of all regular quasimultipliers on $\Uu$ is a pseudo-Banach algebra in the sense of \cite{ad}.

In section 6, we associate in a natural way  to a quasimultiplier $S$ on $\iit$ a weakly densely defined operator $\tilde S$ on $X_1.$ If $\cap_{t>0}Ker(T(t))=\{0\}$ then $\tilde S$ has an extension $\overline{\tilde S}$ to a subspace of $X_1$ which is a closed operator. In the general case the domain of $\tilde S$ can be considered as a subspace of $X_0:=X_1/\cap_{t>0}Ker(T(t))$ and $\tilde S$ has an extension $\overline{\tilde S}$ to a subspace of $X_0$ which is a closed operator.

In section 7 we use a variant of a construction introduced in \cite{e_2} to define the "normalized Arveson ideal" $\jjt.$ Set $\rho_{\bf T}:=\lim_{t\to +\infty}\Vert T(t)\Vert^{1/t}.$ If $\lambda > log(\rho_{\bf T}),$ then the family $\{e^{-\lambda t}T(t)\}_{t>0}\}$ is pseudobounded. Set $\Vert u\Vert_{\lambda} =\sup_{t>0}\Vert e^{-\lambda t}T(t)u\Vert\in [0, +\infty]$ for $u \in \iit,$ set $\lt :=\{ u \in \iit \ | \ \Vert u \Vert_{\lambda} < +\infty\},$ and denote by $\dt$ the closure of $\cup_{t>0}T(t)\iit$ in $(\lt, \Vert .\Vert_{\lambda}).$ Then $(\lt, \Vert .\Vert)_{\lambda}$ is a Banach algebra, and the definition of $\lt$ and $\dt$ and the norm topology on $\lt$ and $\dt$ do not depend on the choice of $\lambda > log(\rho_{\bf T}).$ We can identify the usual multiplier algebra ${\mathcal M}(\dt),$ which is a Banach algebra with respect to the operator norm $\Vert S\Vert_{op, \lambda},$ to a subalgebra of $\qm_r(\iit).$
The closure  $\jjt$ of $\iit$ in $({\mathcal M}(\dt), \Vert S\Vert_{op, \lambda})$ is called the normalized Arveson ideal of the semigroup ${\bf T}.$ The sequence $(\phi_{\bf T}(f_n*\delta_{\epsilon_n}))_{n\ge 1}$ is a sequential bounded approximate identity for $\jjt$ for every Dirac sequence $(f_n)_{n\ge 1}$ and every sequence $(\epsilon_n)_{n\ge 1}$ of positive real numbers which converges to $0,$ and the map $S \to S_{|_{\dt}}$ is an isometric isomorphism from $({\mathcal M}(\jjt), \Vert . \Vert_{{\mathcal M}(\jjt)})$ onto $({\mathcal M}(\dt), \Vert .\Vert_{\lambda, op}).$ The algebras $\dt, \iit$ and $\jjt$ are similar in the sense of \cite{e2}, and an explicit pseudobounded isomorphism between the quasimultiplier algebras $\qm(\jjt)$ and $\qm(\dt)$ or $\qm(\iit)$ (resp. between the algebras $\qm_r(\jjt)$ and $\qm_r(\dt)$ or $\qm_r(\iit)$ ) is given by the map $S_{u/v} \to S_{au/av}$ where $a \in \Omega(\dt).$ We have $\Omega(\jjt)\cap\dt \subset \Omega(\dt),$ and the fact that $\Omega(\jjt)\cap \dt \neq \emptyset$  follows from a construction based on inverse Laplace transforms of outer functions on half-planes given by P. Koosis and the author in section 6 of \cite{e1}.

In section 8 we  consider the generator $A_{\bf T}$ of the semigroup as a quasimultiplier on the algebra $\iit$ or $\jjt.$ To define the infinitesimal generator as a quasimultiplier on $\iit$ we use the formula

$$A_{\bf T, \iit}=S_ {-{\phi_{\bf T}(f_0')/\phi_{\bf T}(f_0)}},$$

where $f_0 \in {\mathcal C}^1([0,+\infty))\cap \Omega \left (L^1_{\omega_{\bf T}}\right )$ satisfies $f_0=0,$ $f_0' \in L^1_{\omega_{\bf T}}.$

This definition does not depend on the choice of $f_0$ and the quasimultiplier $A_{\bf T, \iit}$ is a closed operator on the Arveson ideal $\iit.$ If $u \in \iit,$ and if we have $\lim_{t\to 0^+}\Vert {T(t)u-u\over t}-v\Vert=0$ for some $v \in \iit,$ then $u$ belongs to the domain of $A_{\bf T, \iit}$ and $A_{\bf T, \iit}u=v.$ It is also possible to consider the infinitesimal generator of the semigroup as a quasimultiplier on the normalized Arveson ideal $\jjt.$ Set $\tilde \omega_{\bf T} (t)=\Vert T(t)\Vert_{{\mathcal M}(\jjt)},$ and define $\tilde \phi_{\bf T}= {\mathcal M}_{\tilde \omega_{\bf T}} \to {\mathcal M}(\jjt)$ by using the formula

$$\tilde \phi_{\bf T}(\mu):=\int_0^{+\infty}T(t)d\mu(t),$$

where the Bochner integral is computed with respect to the strong operator topology on ${\mathcal M}(\jjt).$ To define the infinitesimal generator as a quasimultiplier on $\jjt$ we use the formula

$$A_{\bf T, \jjt }=S_ {-{\phi_{\bf T}(f_0')/\phi_{\bf T}(f_0)}},$$

where $f_0 \in {\mathcal C}^1([0,+\infty))\cap \Omega \left (L^1_{\tilde \omega_{\bf T}}\right )$ satisfies $f_0=0,$ $f_0' \in L^1_{\omega_{\bf T}}.$

Notice that while the fact that ${\mathcal C}^1([0,+\infty))\cap \Omega \left (L^1_{\omega_{\bf T}}\right )\neq \emptyset$ follows from nontrivial results from section 6 of \cite{e1} it is very easy to find elements $f_0 \in {\mathcal C}^1([0,+\infty))\cap \Omega \left (L^1_{\tilde \omega_{\bf T}}\right )$ satisfying $f_0=0,$ $f_0' \in L^1_{\omega_{\bf T}}:$ one can take for example $f_0=v_{\lambda},$ for some $\lambda > log(\rho_{\bf T}),$ where $v_{\lambda}(t)=te^{-\lambda t}$ for $t\ge 0.$ This definition agrees with the classical definition of the genarator of a strongly continous semigroup bounded in norm near $0:$ if $u \in \jjt,$ then $u$ belongs to the domain of $A_{\bf T, \jjt }$ if and only if there exists $v \in \jjt$ such that 
$\lim_{t\to 0^+}\Vert {T(t)u-u\over t}-v\Vert_{\jjt}=0,$ and in this situation $A_{\bf T, \jjt}u=v.$

Identify the multiplier algebras $\qm(\iit)$ and $\qm(\jjt),$ and set $A_{\bf T}= A_{\bf T, \iit}=A_{\bf T, \jjt}.$ Every character $\chi$ on $\iit$ extends in a unique way to  a character $\tilde \chi$ on $\qm (\iit),$  and there exists a unique complex number $a_{\chi}=\tilde \chi(A_{\bf T})$ such that $\chi(T(t))=e^{ta_{\chi}}$ for $t>0.$ If the Arveson ideal $\iit$  is not a radical algebra, then the family $\sigma_{ar}(S):= \{ \tilde \chi (S)\}_{\chi \in \sia}$ wil be called the Arveson spectrum of a quasimultiplier $s \in \qm(\iit),$ and we will use the convention
$\sigma_{ar}(A_{\bf T})=\emptyset$ if  the algebra $\iit$ is radical. If $\iit$ is not radical, the map $\chi \to \tilde \chi (A_{\bf T})$ defines a homeomorphism from $\sia$ onto $\sigma_{ar}(A_{\bf T}),$ and we have, for $f \in {\mathcal M}_{\omega_{\bf T}},$

$$\tilde \chi (\phi_{\bf T}(f))={\mathcal L}(f)(-\tilde \chi(A_{\bf T})).$$

Similarly, we have $\tilde \chi (\tilde \phi_{\bf T}(f))={\mathcal L}(f)(-\tilde \chi(A_{\bf T}))$ for $f\in {\mathcal M}_{\tilde \omega_{\bf T}},$ see section 9.

In section 10 we introduce the Arveson resolvent set $Res_{ar}(A_{\bf T}):=\{\C \setminus \sigma_{ar}(A_{\bf T})\},$ and observe that $A_{\bf T}-\mu I$ has an inverse $(A_{\bf T}-\mu I)^{-1}$ in $\qm(\iit)$ for $\mu \in Res_{ar}(A_{\bf T}).$ More precisely this quasimultiplier $(A_{\bf T}-\mu I)^{-1})$ belongs to $\jjt \subset \qm_r(\iit),$ and the map $\mu \to (A_{\bf T}-\mu I)^{-1})$ is an holomorphic map from $Res_{ar}(A_{\bf T})$ into $\jjt.$ If $Re(\lambda)> -log(\rho_{\bf T}),$ then we have the usual resolvent formula

$$(A-\lambda I)^{-1}=-\int_0^{+\infty}e^{-\lambda t}T(t)dt \in \jjt,$$

where the Bochner integral is computed with respect to the strong operator topology on ${\mathcal M}(\jjt).$

We can also define in this case $(A-\lambda I)^{-1}$ by using the formula

$$(A-\lambda I)^{-1}v=-\int_0^{+\infty}e^{-\lambda t}T(t)vdt \ \ (v \in \dt),$$

which defines a quasimultiplier on $\iit$ if we apply it with $v \in \dt \cap\Omega(\iit).$

In section 11 we introduce a holomorphic functional calculus. For $\alpha \in \R,$ set $\Pi^+_{\alpha}:=\{ z \in \C \ | \ Re(z)<\alpha\}.$ Let $H^1(\Pi_{\alpha}^+)$ be the usual Hardy space consisting in those  functions $F$ holomorphic on $\Pi^+_{\alpha}$ such that $\Vert F\Vert_1:= \sup_{\beta >\alpha}{1\over 2\pi}\int_{-\infty}^{+\infty}\vert F(\beta +iy)\vert dy <+\infty.$ Functions $F\in  H^1(\Pi_{\alpha}^+)$ admit a.e. a nontangential limit  $F^*(\alpha +iy)$ on $\alpha +i\R,$ and $\Vert F\Vert_1:={1\over 2\pi}\int_{-\infty}^{+\infty}\vert F^*(\alpha +iy)\vert dy. $ The restriction to $\Pi_{\beta}^+$ of a function $F\in  H^1(\Pi_{\alpha}^+)$ is bounded on $\Pi_{\beta}^+$ for $\beta >\alpha,$ and so $\cup_{\alpha < -log(\rho_{\bf T})}H^1(\Pi_{\alpha}^+)$ is an algebra.

 Set, for $F \in H^1(\Pi_{\alpha}^+),$ $\alpha \in (-\infty, - log(\rho_{\bf T})),$

\begin{center}$F(-A_{\bf T})=-{1\over 2\pi}\int_{-\infty}^{+\infty}F^*(\alpha +iy)(A_{\bf T} + (\alpha +iy)I)^{-1}dy \in \jjt \subset {\mathcal QM}_r(\iit).$\end{center}

Then $FG(-A_{\bf T})=F(-A_{\bf T})G(-A_{\bf T}), F_{|_{\Pi_{\beta}^+}}(-A_{\bf T})$ for $\beta \in (\alpha, -log(\rho_{\bf T})),$ and $F(-A_{\bf T})=\phi_{\bf T}({\mathcal L}^{-1}(F))$ for $
F \in \cup_{\alpha < log(\rho_{\bf T})} H^1(\Pi_{\alpha}^+),$ $G \in \cup_{\alpha < log(\rho_{\bf T})}H^1(\Pi_{\alpha}^+).$

After observing in section 11 that $G(-A_{\bf T})\in \Omega(\jjt)$ if $G \in H^1(\Pi^+_{\alpha})$ is outer (this can be deduced from Nyman's characterization of dense ideals of the convolution algebra $L^1(\R^+),$ but we propose a direct proof), we define $F(-A_{\bf T})$ for $F \in \cup_{\alpha < log(\rho_{\bf T})}H^{\infty}(\Pi_{\alpha}^+)$ by using the formula

$$F(-A_{\bf T})=S_{FH(-A_{\bf T})/H(-A_{\bf T})}\in \qm_r(\jjt)=\qm_r(\iit),$$

where $F\in H^{\infty}(\Pi_{\alpha}^+), FH\in H^{1}(\Pi_{\alpha}^+),$ and where $H\in H^{1}(\Pi_{\alpha}^+)$ is outer on $\Pi_{\alpha}^+.$

Again, this definition does not depend on the choice of $H.$ When $F={\mathcal L}(f),$ with $f\in \cup_{\lambda >log(\rho_{\bf T})}L^1(\R^+,e^{\lambda t}),$ then $F(-A_{\bf T})=\tilde \phi_{\bf T}(f)\in \jjt.$ If $F(z)=e^{-tz},$ with $t>0,$ then $F(-A_{\bf T})=T(t).$ A similar formula allows to define $F(-A_{\bf T})\in \qm(\iit)=\qm(\jjt),$ if $F\in \cup_{\alpha < log(\rho_{\bf T})} {\mathcal N}^+(\Pi_{\alpha}^+),$ where ${\mathcal N}^+(\Pi_{\alpha}^+)$ denotes the Smirnov class on $\Pi^+_{\alpha},$ which is the class of all holomorphic functions $F$ on $\Pi_{\alpha}^+$ such that $FH \in H^{\infty}(\Pi_{\alpha}^+)$ for some outer function $F\in H^{\infty}(\Pi_{\alpha}^+).$ We have again $FG(-A_{\bf T})=F(-A_{\bf T})G(-A_{\bf T}),$ and, indeed, $F(-A_{\bf T})=A_{\bf T}$ if $F(z)=-z.$

I. Chalendar, J. Partington and the author obtained in \cite{cep1} lower estimates for small values of $\epsilon$ of the norm of operators of the form $F(-\epsilon A_{\bf T}),$ where $F$  is the Laplace transform of a real-valued measure on $[a,b],$ with $0<a<b<+\infty,$ and where ${\bf T}=(T(t))_{t>0}$ is a norm-continuous semigroup, see also \cite{em},\cite{e3}, \cite{bcep}. The motivation for this paper was to open the road to extensions of these lower estimates to much more general situations.

The idea of considering the infinitesimal generator of a semigroup as a quasimultipier on a suitable Banach algebra associated to the semigroup is due to Gal\'e and Miana, who used this approach in \cite{gm} for the infinitesimal generator of some one-parameter groups of regular quasimultipliers, see also \cite{ga}. Chalendar, Partington and the author also used this approach for the infinitesimal generator on analytic semigroups. We deal here with a more general situation, and the author hopes that this approach centered on quasimultipliers and regular quasimultipliers on the Arveson ideal $\iit$ will be helpful to mathematicians who use holomorphic functional calculus associated to semigroups.

\section {Arveson weakly continuous one-parameter semigroups}

     Let $X=(X,\Vert . \Vert)$ be a Banach space. We denote by ${\mathcal B}(X)$ the 
Banach algebra of bounded linear operators $R:X\to X$ with composition, 
we denote by ${\mathcal{GL}}(X)$ the group of invertible elements of ${\mathcal B}(X)$ and we denote by  $I=I_X$  the identity map on $X.$ We also denote by $\Vert . \Vert$ the operator norm on  ${\mathcal B}(X)$ associated to the given norm on $X,$ and we denote by $\rho( R )$ the spectral radius of $R \in {\mathcal B}(X).$ If $Y$ is a subspace of the dual space $X'$ of $X$ we will denote by $\sigma(X,Y)$ the weak topology on $X$ associated to $Y.$
     
     We will use  the following notion, introduced by Arveson in \cite{ar}.
     
     \begin{defn} \label{dualpairing} Let $X$ be a Banach space, and let $X_*$ be a 
subspace of the dual space $X'.$ We will say that $(X,X_*)$ is a dual pairing if the 
two following conditions are satisfied:
     
     \smallskip
     \begin{enumerate}
     
     \item $\Vert x \Vert=\sup \{\vert  \langle x,l\rangle \vert \  :  \ {l \in X_*}, \Vert l \Vert \le 1\}$ for every $x \in X.$
     
     \smallskip
     
     \item The $\sigma(X,X_*)$-closed convex hull of every $\sigma(X,X_*)$-compact subset of $X$ is $\sigma(X,X_*)$-compact.
     
     \smallskip
     \end{enumerate}
     \end{defn}
     \noindent 
     For example, $(X,X')$ is a dual pairing. Also if $X=Y'$ for some Banach space $Y,$ and if we identify $Y$ to a subspace of $X'=(Y')'$ in the obvious way, then $(X,Y)$ is a dual pairing, 
see \cite{ar}. Notice that condition 1 means that if we set $\tilde x(l)=\langle x,l\rangle $ for $x \in X, l \in X_*$ 
then the map $x \mapsto \tilde x$ is an isometry from $X$ into the dual space $(X_*)'.$ This means that the spaces $X$ and $X_*$ are reciprocal in the sense of definition 2 of \cite{f0}. An operator $T \in {\mathcal B}(X)$ is said to be weakly continuous with respect to a dual pairing $(X,X_*)$ if $l\circ T \in X_*$ for every $l\in X_*.$ 

     Let $S$ be a locally compact space, and let $(X,X_*)$ be a dual pairing. 
A map $u: S \to X$ is said to be {\it weakly continuous with respect to}  $(X,X_*)$ if the map 
$s \mapsto \langle u(s),l\rangle $ is continuous on $S$ for every $l \in X_*,$ and we will often just say that 
$u$ is weakly continuous when no confusion may occur. In this situation it follows from the 
Banach-Steinhaus theorem and from condition 1 of definition \ref{dualpairing} that we have, for every compact subset $K$ of $S,$     
     \begin{equation} \sup_{s \in K}\Vert u(s)\Vert=\sup_{s \in K}\Vert 
\widetilde{ u (s)}\Vert
 = \sup_{s\in K}\sup_{l\in X_*,\Vert l \Vert \leq 1} \vert \langle  u(s),l\rangle \vert <
+\infty .\end{equation}

Since $\Vert u(s) \Vert =\sup_{l \in X_*}\vert\langle  u(s),l\rangle \vert$ for $s \in S,$ the function 
$\omega_u:s \mapsto \Vert u(s)\Vert $ is lower semicontinuous on $S,$ 
which allows to compute the upper integral 
\begin{equation}\int^*_S \Vert u(s) \Vert d| \mu|(s):=\sup_{\stackrel{f \in {\mathcal C}^+_c(S)}{_{ f \le \omega_u}}}
\int_S f(s) d| \mu | (s)\in [0,+\infty]
\end{equation} for every regular measure $\mu$ on $S,$ where ${\mathcal C}^+_c(S)$ denotes the space of all nonnegative compactly supported continuous functions on $S.$
     The following proposition is an immediate generalization of proposition 1.2 of \cite{ar}. The details of the proof are given in \cite{ef}.
     
     \begin{prop}  Let $S$ be a locally compact space, let $(X,X_*)$ be a dual pairing, 
and let  $u: S \to X$ be a weakly continuous map. Set $\omega_u(s) =\Vert u(s)\Vert$ for
 $s \in S,$ and denote by ${\mathcal M}_{\omega_{u}}(S)$ the set of all regular measures 
$\mu$ on $S$ 
such that $\Vert \mu \Vert_{\omega_u}:=\int_S^*\Vert u(s)\Vert d|\mu|(s)<+\infty.$ 
Then for every $\mu \in {\mathcal M}_{\omega_{u}}(S)$ there exists $x\in X$ satisfying

\begin{equation} \langle x,l\rangle =\int_S\langle u(s),l\rangle d\mu(s) \ \ \ \ (l \in X_*).\end{equation}
\end{prop}
Proof: Since $\int_S\vert \langle u(s),l\rangle \vert d|\mu|(s) \le \Vert l\Vert 
\Vert \mu \Vert_{\omega_u}<+\infty$ for every $l \in X_*,$ the formula $f_{\mu}(l):=\int_S\langle u(s),l\rangle d\mu(s)$ for  $l \in X_*$  defines an element $f_\mu \in (X_*)',$ and we have to show that $f_{\mu}=\tilde x$ for some $x\in X.$ 
It follows from condition 1 of definition \ref{dualpairing} that we have
$$\Vert f_{\mu}\Vert=\sup_{l\in X_*,\Vert l\Vert \le 1}\left | \int_S\langle u(s),l\rangle d\mu(s) 
\right | \le \Vert \mu \Vert_{\omega_u}.$$
Denote by ${\mathcal M}_c(S)$ the space of all regular measures on $S$ supported by some 
compact subset of $S.$ It follows from (2)  that 
${\mathcal M}_c(S)\subset {\mathcal M}_{\omega_u}(S),$ and the fact that property 
(2.2) holds for every $\mu \in {\mathcal M}_c(S)$ follows directly from \cite{ar}, 
proposition 2.1. Set $\tilde X:=\{ \tilde x : x \in X\}.$ It follows from condition 1 
of definition \ref{dualpairing} that $\tilde X$ is closed in $(X_*)'.$ 
Let $\mu \in {\mathcal M}_{\omega_u}(S).$ There exists a sequence $(\mu_n)_{n\ge 1}$ 
of elements of ${\mathcal  M}_c(S)$ such that 
$\lim_{n\to +\infty}\Vert \mu -\mu_n\Vert_{\omega_u}=0.$ Hence 
$\lim_{n\to +\infty}\Vert f_{\mu}-f_{\mu_n}\Vert=0,$ and $f_\mu \in \tilde X.$ \hfill $\square$

\vskip 2 mm 

When the conditions of proposition 2.2 are satisfied, we will use the notation
\begin{equation} x = \int_Su(s)d\mu(s),
\end{equation}
where the integral is a Pettis integral, that is a Bochner integral computed with respect to the $\sigma(X,X_*)$ topology, 
which defines an element of $X$ since $(X,X_*)$ is a dual pairing.

Notice that since $\omega_u$ is lower semicontinuous, it follows from the theory of integration 
on locally compact spaces, see \cite{bo}, chapter 4, that ${\mathcal M}_{\omega_u}(S)$ is the 
space of regular measures $\mu$ such that $\omega_u$ is integrable with respect to the total 
variation $\vert \mu \vert$ of $\mu,$ and it follows from \cite{bo}, proposition 1, that we have 
for $\mu \in {\mathcal M}_{\omega_u}(S)$
$$\Vert \mu\Vert _{\omega_u}= \int _S \Vert u(s)\Vert d| \mu | (s) 
=\int^*_S\Vert u(s) \Vert d|\mu|(s)= \sup_{K\subset S, K { compact}}\int^*_K\Vert u(s)\Vert d | \mu|(s).$$

Let $X$ be a Banach space. A map ${\bf T}: t \to T(t)$ from $(0,+\infty)$ into ${\mathcal B}(X)$ is called a (one-parameter) semigroup if $T(s+t) =T(s)T(t)$ for $s>0,t>0.$ We will often write ${\bf T}=(T(t))_{t>0}.$ A semigroup will be said to be normalized if $\cup_{t>0}T(t)(X)$ is dense in $X$.

We now introduce the following definition, which is an obvious extension to one-parameter semigroups of the notion of weakly continuous group representations associated to a dual pairing, due to Arveson \cite{ar}.

\begin{defn} Let $(X,X_*)$ be a dual pairing, and let ${\mathcal B}_w(X)$ the the space of weakly continuous elements of $B(X).$ A semigroup ${\bf T}=(T(t))_{t>0} \subset {\mathcal B_w}(X)$ is said to be weakly continuous with respect to the dual pairing $(X,X_*)$ if the function $t \to <T(t)x,l>$ is continuous on $(0,+\infty)$ for every $x\in X$ and every $l \in X_*.$
\end{defn}

Let $(X,X_*)$ be a dual pairing, and let ${\bf T}=(T(t))_{t>0}\subset {\mathcal B}_w(X)$ be a weakly continuous 
semigroup. We have, for $t>0,$
$$\Vert T(t) \Vert =\sup\{\vert  \langle T(t)x,l\rangle \vert : x \in X, \Vert x \Vert \le 1, 
l \in X_*, \Vert l \Vert \le 1\}$$
and so the weight $\omega_{\bf T}: t \mapsto \Vert T(t)\Vert$ is lower semicontinuous on $(0,+\infty).$ 
Let $K\subset (0,+\infty)$ be compact. Since $\sup_{t \in K}\Vert T(t)x \Vert <+\infty$ for every 
$x \in X,$ it follows again from the Banach-Steinhaus theorem that 
$\sup_{t \in K}\Vert T(t)\Vert <+\infty.$ A standard argument shows then that $lim_{t\to +\infty}\Vert T(t)\Vert^{1\over t}= lim_{n\to +\infty}\Vert T(sn)\Vert^{1\over sn}=inf_{n\ge 1}\Vert T(sn)\Vert^{1\over sn}$ for every $s>0.$ In this situation we can define the weighted space 
${\mathcal M}_{\omega_{\bf T}}={\mathcal M}_{\omega_{\bf T}}[(0,+\infty)]$ consisting of all regular measures $\mu$ on $(0,+\infty)$ 
such that the upper integral $\int_0^{+\infty}\Vert T(t)\Vert d\vert \mu \vert(t)$ is finite. 
Since $\omega_{\bf T}(t_1+t_2)\le \omega_{\bf T}(t_1)\omega_{\bf T}(t_2)$ for $t_1>0, t_2 >0,$ 
we see that $\mu$ and $\nu$ are convolable for  $\mu,\nu \in {\mathcal M}_{\omega_{\bf T}}$ and that $({\mathcal M}_{\omega_{\bf T}}, \Vert . \Vert_{\omega_{\bf T}})$ is a Banach algebra with respect to convolution which contains the convolution algebra ${\mathcal M}_c[(0,+\infty)]$ of compactly supported regular measures on $(0,+\infty)$ as a dense subalgebra.

Denote by $L^1_{\omega_{\bf T}}=L_{\omega_{\bf T}}^1(\R^+)$ the convolution algebra of all Lebesgue-measurable (classes of) functions $f$ on $[0,+\infty)$ such that $f\omega_{\bf T}$ is integrable with respect to the Lebesgue  measure on $[0,+\infty),$ identified to the space of all measures $\mu \in {\mathcal M}_{\omega_{\bf T}}$ which are absolutely continuous with respect to Lebesgue measure. Also denote by ${\mathcal C}^{\infty}_c[(0,+\infty)]$ the space of all infinitely differentiable functions $f \in {\mathcal C}_c[(0,+\infty)].$
We will say as usual that a sequence $(f_n)_{n\ge 1}$ of elements of $L^1(\R^+)$ is a {\it Dirac sequence} if $f_n(t)\ge 0$ a.e. for $t\ge 0$, $\int_0^{+\infty}f_n(t)dt=1$ and if there exists a sequence $(a_n)_{n\ge 1}$ of positive real numbers which converges to $0$ such that $f_n(t)=0$ a.e. on $[a_n,+\infty)$ for $n\ge 1$.

The following result is an easy analog of \cite{ar}, proposition 1.4.

\begin{prop} Let $(X,X_*)$ be a dual pairing,
 and let ${\bf T}=(T(t)_{t >0} \subset {\mathcal B}_w(X)$ be a weakly continuous semigroup.
The Pettis integral
\begin{equation} <\phi_{\bf T}(\mu)x,l>=\int_0^{+\infty}<T(t)x,>d\mu(t) \ \  (l\in X_*)\end{equation}
defines for every $x \in X$ and every $\mu \in {\mathcal M}_{\omega_{\bf T}}$ an element of 
$X,$ $\phi_{\bf T}(\mu)\in {\mathcal B}(X)$ for every $\mu \in {\mathcal M}_{\omega_{\bf T}},$ 
and $\phi_{\bf T}: \mu \mapsto \phi_{\bf T}(\mu)$ is a norm-decreasing unital algebra homomorphism 
from the convolution algebra ${\mathcal M}_{\omega_{\bf T}}$ into ${\mathcal B}(X).$

 Moreover if  $(f_n)_{n\ge 1} $ is a Dirac sequence,  then we have, for  $t>0$

$$(i) lim_{n\to +\infty}\Vert \phi_{\bf T}(f_n*\delta_t)\phi_{\bf T}(f)-T(t) \phi_{\bf T}(f)\Vert =0  \ \forall f \in \LL,$$

$$(ii) lim_{n \to +\infty} \langle  \phi_{\bf T}(f_n*\delta_t)x-T(t)x,l\rangle =0 \ \forall x \in X, \forall l \in X_*.$$
\end{prop}

Proof: Since $\Vert T(t)x \Vert \le \Vert T(t) \Vert \Vert x \Vert,$ the fact that formula 
5  defines an element of $X$ for $x \in X$ and 
$\mu \in {\mathcal M}_{\omega_{\bf T}}$ follows directly from proposition 2.2. We have
\begin{eqnarray*}
\Vert \phi_{\bf T}(\mu)x \Vert &=&\sup_{l \in X_*, \Vert l \Vert \le 1}\left | 
\int_0^{+\infty} \langle T(t)x,l\rangle d\mu(t)\right | \\
&\le& \sup_{l \in X_*, \Vert l \Vert \le 1}\int_0^{+\infty}\left | \langle  T(t)x,l\rangle \right  |d | \mu |(t)
\le \Vert x \Vert \Vert \mu \Vert_{\omega_{\bf T}},
\end{eqnarray*} 
and so $\phi_{\bf T}(\mu) \in {\mathcal B}(X)$ for $\mu \in {\mathcal M}_{\omega_{\bf T}}$ 
and $\Vert \phi_{\bf T}(\mu)\Vert \le \Vert \mu \Vert_{\omega_{\bf T}}.$    
As observed in \cite{ar}, a routine application of Fubini's theorem shows that $\phi_{\bf T}(\mu*\nu)=\phi_{\bf T}(\mu) \phi_{\bf T}(\nu)$ for $\mu, \nu \in {\mathcal M}_c((0,+\infty)).$ Since ${\mathcal M}_c((0,+\infty))$ is dense in ${\mathcal M}_{\omega_{\bf T}},$ this shows that $\phi_{\bf T}$ is an algebra homomorphism. 

The last assertions pertain to folklore, but we give the details for the sake of completeness. 
Let $n\ge 1,$ and let $(f_n)_{n\ge 1}$ be a Dirac sequence. We can assume without loss of generality that $f_n(t)=0$ a.e. on $[1,+\infty)$ for $n \ge 1.$ The sequence $(f_n)_{n\ge 1}$ is a sequential bounded approximate identity for $L^1(\R^+).$ So if $f \in \LL,$ and if $supp(f)\subset [a, b],$ with $0<a<b<+\infty,$ we have

$$lim sup_{n\to +\infty}\Vert f_n*\delta_t*f -\delta _t *f\Vert_{\LL}\le lim Msup_{n\to +\infty}\Vert f_n*\delta_t*f -\delta _t *f\Vert_{L^1(\R^+)}=0,$$

where $M=sup_{a+t \le s \le b+1+t}\Vert \omega_{\bf T}(s)\Vert<+\infty.$

Since $\Vert f_n*\delta_t\Vert_{\LL}\le sup_{t \le s \le t+1}\Vert T(s)\Vert\int_{0}^1f_n(s)ds=sup_{t \le s \le t+1}\Vert T(s)\Vert$ for $n\ge 1,$ a standard density argument shows that $lim sup_{n\to +\infty}\Vert f_n*\delta_t*f -\delta _t *f\Vert_{\LL}=0,$ and (i) holds.

We have, for $x \in X, l \in X_*,t>0$
\begin{eqnarray*}
 \left | \langle \phi_{\bf T}(f_n*\delta_t)x -T(t)x,l\rangle \right | &=& \left | \int_0^{a_n}\langle T(s+t)x ,l\rangle f_n(s)ds -
     \langle T(t) x,l\rangle \int_0^{a_n}f_n(s)ds\right | \\
 &\le& \sup_{0\le t \le {a_n}}\left | \langle  T(s+t)x-T(t)x,l\rangle \right |\int_0^{a_n}f_n(t)dt)\\
&=& \sup_{0\le t \le {a_n}}\left | \langle  T(s+t)x-T(t)x,l\rangle \right |.
\end{eqnarray*} Hence $lim_{n \to +\infty}  \left | \langle \phi_{\bf T}(f_n*\delta_t)x -T(t)x,l\rangle \right | =0,$ since the semigroup is weakly continuous with respect to the dual pairing $(X,X_*).$
$\square$

In the following $(X,X_*)$ will denote a dual pairing and ${\bf T}=(T(t))_{t>0}\subset {\mathcal B}_w(X)$ will denote a  semigroup of bounded weakly continuous operators on $X$ which is weakly continuous with respect to this dual pairing $(X,X_*).$ The weak topology $\sigma(X,X_*)$ associated to the dual pairing will be called the weak topology on $X.$ 

We will use  the following notations, for $x\in X,$ $\mu \in \M,$

\begin{equation} \phi_{\bf T}(\mu)x= \int_0^{+\infty}T(t)xd\mu(t), \phi_{\bf T}(\mu)= \int_0^{+\infty}T(t)d\mu(t).\end{equation}

In particular we will write, for $f \in \LL,$

\begin{equation} \int_{0}^{+\infty}f(t)T(t)dt=\phi_{\bf T}(f).\end{equation}

\section{The Arveson ideal}

Denote by  $X_1=\left [ \cup_{t>0}T(t)(X)\right ]^{-}$ the closure of the space $\cup_{t>0}T(t)(X)$ in $(X,\Vert .\Vert).$ Since ${\mathcal M}_c((0,+\infty))$ is dense in ${\mathcal M}_{\bf T},$ we have $\phi_{\bf T}(\mu)x\in X_1$ for $\mu \in {\mathcal M}_{\omega_{\bf T}},$  $x\in X.$ So we may consider
$\phi_{\bf T}$ as an element of ${\mathcal B}(X_1)$ for every $\mu \in {\mathcal M}_{\omega_{\bf T}},$ and $\Vert \phi_{\bf T}(\mu)\Vert _{{\mathcal B}(X)}\ge \Vert \phi_{\bf T}(\mu)\Vert _{{\mathcal B}(X_1)}$ for every $\mu \in {\mathcal M}_{\bf T}.$ The trivial example obtained by setting $T(t)=P$ for $t>0,$ where $P:X \to Im(P)$ is a projection satisfying $\Vert P\Vert_{{\mathcal B}(X)}>1,$ shows that equality does not hold in general.

\begin{defn} We will denote by ${\mathcal A}_{\bf T}$ the closed subalgebra of ${\mathcal B}(X_1)$ generated by $\phi_{T}({\mathcal M}_{\omega_{\bf T}}),$ and by $\iit$ the closed subalgebra of 
${\mathcal B}(X_1)$ generated by $\phi_{T}({\LL}),$ so that $\iit$ is a closed  ideal of ${\mathcal A}_{\bf T},$ called the Arveson ideal of the semigroup ${\bf T}.$
\end{defn}

The fact that $\iit $ is an ideal of  ${\mathcal A}_{\bf T}$ follows from the fact that $\LL$ is an ideal of the convolution algebra ${\mathcal M}_{\omega_{\bf T}},$ so that  $\phi_{T}({\LL})$
is an ideal of $\phi_{T}({\mathcal M}_{\omega_{\bf T}}).$

\begin{prop} (i) The map $t \to T(t)u$ is norm-continuous on $(0,+\infty)$ for every $u \in \iit.$

\smallskip

(ii) The ideal $\iit$ equals the whole algebra ${\mathcal A}_{\bf T}$ if and only if the semigroup $(T(t))_{t>0}$ is continuous with respect to the norm of ${\mathcal B}(X).$

\smallskip

 (iii) If the $\sigma(X_*,X)$-closed convex hull of every $\sigma(X_*,X)$-compact subset of $X$ is $\sigma(X_*,X)$-compact, then all elements of $\at$ are weakly continuous with respect to $X_*.$

\end{prop}

Proof: The map $t \to \delta_t*f$ is continuous on $(0,+\infty)$ for every $f \in \LL.$ Hence the map $t \to T(t)\phi_{\bf T}(f)$ is continuous on $(0,+\infty)$ for every $f \in \LL.$
Since $sup_{a\le t \le b}\Vert T(t)\Vert<+\infty$ for $0<a<b<+\infty,$ this shows that (i) holds.

Now assume that the semigroup is continuous with respect to the norm on ${\mathcal B}(X).$ A well-known argument, see for example  proposition 6.1 in \cite{e1}, shows then that $T(t)\in \iit$ for every $t>0.$ Hence $\phi_{\bf T}(\mu)\in \iit$ for every $\mu \in {\mathcal M}_c((0,+\infty)).$ Since ${\mathcal M}_c((0,+\infty))$ is dense in $\M,$ this shows that ${\mathcal A}_{\bf T}=\iit.$ 

Conversely it follows immediately from (i) that if ${\iit}= {\mathcal A}_{\bf T}$ then the semigroup $(T(t))_{t>0}$ is continuous with respect to the norm of ${\mathcal B}(\tilde X).$ 

Now assume that the $\sigma(X_*,X)$-closed convex hull of every $\sigma(X_*,X)$-compact subset of $X$ is $\sigma(X_*,X)$-compact. Then $(X_*,X)$ is a dual pairing, and we see as in the proof of the second assertion of proposition 1.4 of \cite{ar} that $l\circ \phi_{\bf T}(\mu)\in X_*$ for  $\mu \in {\mathcal M}_{\omega_{\bf T}}, l \in X_*.$ Hence $l\circ u \in X_*$ for $u \in \at, l \in X_*.$
$\square$

 \begin{prop} Let $u\in {\at}.$ If $uv =0$ for every $v\in \iit,$ then $u=0.$

\end{prop}

Proof: Let $(f_n)_{n\ge 1}$ be a Dirac sequence, and let $\epsilon_n$ be a sequence of positive real numbers such that $\lim_{n\to +\infty}\epsilon_n=0.$ We have, for $x \in X, l \in X_*, t>0,$

$$<uT(t)x,l>=<T(t)ux,l>=lim_{n\to +\infty}<\phi_{T}(f_n*\delta_{\epsilon_n}*\delta_t)ux,l>=$$ $$=\lim_{n\to +\infty}<u\phi_{\bf T}(f_n*\delta_{\epsilon_n})T(t)x,l>=0.$$

Hence $uT(t)x=0$ for every $x\in X,$ and so $u=0$ since $\cup_{t>0}T(t)X$ is dense in $ X_1.$ $\square$

Recall that a sequence $(e_n)_{n\ge 1}$ of elements of a commutative Banach algebra ${\mathcal A}$ is said to be a sequential bounded approximate identity for ${\mathcal A}$ if
$sup_{n\ge 1}\Vert e_n\Vert<+\infty$ and if $lim_{n\to +\infty}\Vert ae_n-a\Vert=0$ for every $a\in {\mathcal A}.$

\begin{prop}


  If $lim sup_{t\to 0^+}\Vert T(t) \Vert <+\infty,$ then the sequence $\phi_{\bf T}(f_n)_{n\ge 1}$ is a bounded sequential approximate identity for $\iit$  for every Dirac sequence $(f_n)_{n\ge 1}.$
 
 \end{prop}
 
 Proof: Assume that $lim sup_{t\to 0^+}\Vert T(t) \Vert <+\infty,$ and set $M=sup_{0<t\le 1}\Vert T(t) \Vert.$ Then $\Vert \phi_{\bf T}(f_n)\Vert\le M$ when $n$ is sufficiently large. Since $(f_n)_{n\ge 1}$ is a sequential bounded approximate identity for $L^1(\R^+),$ and since $sup_{a\le t \le b}\Vert T(t)\Vert<+\infty$ for $0<a<b<+\infty,$ we see that $lim_{n\to +\infty}\Vert f_n*f-f\Vert _{\LL}=0$, so that  $lim_{n\to +\infty}\Vert \phi_{\bf T}(f_n)*\phi_{\bf T}(f) -\phi_{\bf T}(f) \Vert =0$ for $f\in {\mathcal C}_c((0,+\infty)).$ Since the sequence $(\phi_{\bf T}(f_n))_{n\ge 1}$ is bounded, this shows that 
 $lim_{n\to +\infty}\Vert \phi_{\bf T}(f_n)*u -u\Vert=0$ for $u\in \iit.$ $\square$

\section{The generator as a weakly densely defined operator}

We consider as above the operators $(T(t))_{t>0}$ as operators on the Banach space $X_1:=\left [ \cup_{t>0}T(t)X\right ]^-.$ We set

$$\tilde X:=\{ x\in X_1 \ | \ \lim_{t\to 0^+} <T(t)x -x,l>=0 \ \forall l \in X_*\}.$$ We will use the convention $T(0)=I_{X_1},$ the identity map on $X_1.$

Clearly,  $\cup_{t>0}T(t)(X)\subset \tilde X,$ and $[\cap_{t>0}Ker(T(t)]\cap \tilde X=\{0\}.$ Notice that if  $x \in \tilde X,$ it follows from the Banach-Steinhaus theorem that lim sup$_{t\to 0^+}\Vert T(t)x\Vert=$lim sup$_{t\to 0^+}\Vert T(t)x\Vert_{(X_*)'}<+\infty.$ Also if  lim sup$_{t>0}\Vert T(t)\Vert <+\infty,$ then $\tilde X$ is closed, so that $\tilde X=X_1.$

The following definition is a variant of a definition introduced by Feller  in \cite{f0} for semigroups satisfying some weak measurability conditions.

\begin{defn} Denote by  ${\mathcal D}_{{\bf T},op}$ the space of all $x \in X_1$ such that there exists $y\in \tilde X$ satisfying 

\begin{equation}\lim_{t\to 0^+}<{T(t)x-x\over t}-y,l>=0 \ \  \forall l\in X_*.\end{equation}

In this situation we set $A_{{\bf T},op}x=y.$

\end{defn}

Notice that  if $x\in {\mathcal D}_{{\bf T},op},$ then $A_{{\bf T},op}x\in X_1.$ Also  $l\circ T(t)\in X_*$ for every $t>0$ and every $l\in X_*,$ and so when  condition (8) is satisfied by $x$ and $y$ we have $\lim_{t\to 0^+}<{T(t+s)x-T(s)x\over t}-y,l>=0$ for every $l \in X_*$ and every $s>0.$

 If $x \in {\mathcal D}_{{\bf T},op},$ and if $y=A_{{\bf T},op}x,$ then we have, for $l\in X_*,$ $t>0,$

$$<T(t)x-x,l>=\int_0^t<T(s)y,l>dt,$$

and so, according to the notations of formula (4), we have

\begin{equation} T(t)x-x=\int_0^tT(s)yds.\end{equation}

Conversely, if $y \in \tilde X,$ and if $T(t)x-x= \int_0^tT(s)yds$ for $t>0,$ then $x$ and $y$ obviuously satisfy (8) and so $x \in {\mathcal D}_{{\bf T},op},$ and $y=A_{{\bf T},op}x.$


\begin{prop} Set $\rho_{\bf T}=\lim_{t+\infty}\Vert T(t)\Vert^{1\over t},$ and for $Re(\lambda) > log (\rho_{\bf T}),$ $y \in \tilde X,$ define $R_{op}({\bf T}, \lambda)y$ by the formula

\begin{equation} <R_{op}({\bf T}, \lambda)y,l>:=\int_0^{+\infty}e^{-\lambda s}<T(s)y,l>ds \ \ \forall l \in X_*.\end{equation}.

In other terms, according to the notations of formula (6), we have

$$R_{op}({\bf T}, \lambda)y=\int_0^{+\infty}e^{-\lambda s}T(s)yds.$$

Then $R_{op}({\bf T}, \lambda)y \in {\mathcal D}_{{\bf T},op}$ for $y \in \tilde X,$ $R_{op}({\bf T}, \lambda): \tilde X \to {\mathcal D}_{{\bf T},op}$ is one-to-one and onto, and we have

$$\left (\lambda I_{X_1}-A_{{\bf T},op}\right )\circ R_{op}({\bf T}, \lambda)=I_{\tilde X},$$

$$R_{op}({\bf T}, \lambda)\circ \left (\lambda I_{X_1}-A_{{\bf T},op}\right )=I_{{\mathcal D}_{{\bf T},op}}.$$

\end{prop}

Proof: The proof is standard, and we give the details for the sake of completeness. Let $y \in \tilde X.$ Since $\int_0^{+\infty} e^{-Re(\lambda)s}\Vert T(s)y\Vert ds<+\infty,$ $R_{op}({\bf T}, \lambda)y$ is well-defined. Set $x=R_{op}({\bf T}, \lambda)y.$ Since $l*T(t)\in X_*$ for every $l \in X_*$ and every $t>0,$ we have

$$T(t)x=\int_0^{+\infty}e^{-\lambda s}T(t+s)yds=e^{\lambda t}\int_t^{+\infty}T(s)yds,$$ $$ {T(t)x-x\over t}={e^{\lambda t}-1\over t}x -e^{\lambda t}{\int_0^tT(s)yds\over t},$$ where the integral is computed with respect to the weak operator topology on $X$ associated to the dual pairing $(X,X_*).$

We obtain, for $l \in X_*,$

$$\lim_{t\to 0^+}<{T(t)x-x\over t},l>=<\lambda x -y,l>,$$ and so $x \in {\mathcal D}_{{\bf T},op},$ and $A_{{\bf T},op}x= \lambda x -y,$ $y=\lambda x-A_{{\bf T},op}x.$

This shows that $R_{op}({\bf T}, \lambda)(\tilde X)\subset {\mathcal D}_{{\bf T},op},$ and that $\left (\lambda I_{X_1}-A_{{\bf T},op}\right )\circ R_{op}({\bf T}, \lambda)=I_{\tilde X},$ and so $R_{{\bf T},op}(\lambda): \tilde X \to {\mathcal D}_{{\bf T},op}$ is one-to-one.

Conversely let $x\in {\mathcal D}_{{\bf T},op},$ set $u=A_{{\bf T},op}x,$ and set $S(t)=e^{-\lambda t}T(t)$ for $t>0.$ Then $x\in {\mathcal D}_{{\bf S},op},$  and if we set $y=A_{{\bf S},op}x$ we have $y=u-\lambda x.$ We have, for $t>0,$

$$S(t)x-x=\int_0^t S(s)yds=\int_0^te^{-\lambda s}T(s)yds.$$

Since $\lim_{t\to +\infty}\Vert S(t)\Vert=0,$ we obtain

$$x= -\int_0^{+\infty}e^{-\lambda s}T(s)yds=-R_{{\bf T},op}(\lambda)y=\lambda R_{{\bf T},op}(\lambda)x -(R_{{\bf T},op}(\lambda)\circ A_{{\bf T},op})x.$$

Hence $R_{{\bf T},op}(\lambda): \tilde X \to {\mathcal D}_{{\bf T},op}$ is onto, and we have $R_{{\bf T},op}(\lambda)\circ ( \lambda I_{X_1} -A_{{\bf T},op})=I_{{\mathcal D}_{{\bf T},op}}.$ $\square$

\begin{prop} The domain ${\mathcal D}_{{\bf T},op}$ of $A_{{\bf T},op}$ is weakly dense in $X_1.$

\end{prop}

Proof: Let $f \in {\mathcal C}_c^{\infty}([0,+\infty)),$ and denote by $\Vert .\Vert_1$ the usual norm on $L^1(\R^+).$ Then $$\lim_{t\to 0^+}\left \Vert {\delta_t*f-f\over t}-f'\right \Vert_{\omega_{\bf T}}=\lim_{t\to 0^+}\left \Vert {\delta_t*f-f\over t}-f'\right \Vert_1=0,$$ and so $\lim_{t\to 0^+}\left \Vert {T(t)*\phi_{\bf T}(f)-\phi_{\bf T}(f)\over t}-\phi_{\bf T}(f')\right \Vert=0.$ Since $\phi_{\bf T}(f') \subset T(t)\at$ for some $t>0,$ we have $\phi_{\bf T}(f)(X)\subset \tilde X,$ and so $\phi_{\bf T}(f)(X)\subset {\mathcal D}_{{\bf T},op}.$ Now let $(f_n)_{n\ge 1}\subset {\mathcal C}_c^{\infty}([0,+\infty))$ be a Dirac sequence. It follows from proposition 2.4 that we have for $x\in X, l \in X_*,$ $t>0,$

$$\lim_{n\to +\infty}< \phi_{\bf T}(f_n)T(t)x -T(t)x,l>=0.$$

Hence ${\mathcal D}_{{\bf T},op}$ is weakly dense in $\cup_{t>0}T(t)(X),$ and so ${\mathcal D}_{{\bf T},op}$  is weakly dense in $X_1.$ $\square$

\smallskip

 Recall that a partially defined operator $S: {\mathcal D} \to  E$ on a Banach space $E$ is said to be closed if the graph $\{ x, Sx : x \in {\mathcal D}\}$ of $S$ is closed. Notice that if lim sup$_{t\to 0^+}\Vert T(t)\Vert<+\infty,$ then $\tilde X=X_1,$ and so $R({\bf T}, \lambda): \tilde X \to {\mathcal D}_{{\bf T},op}\subset \tilde X$ is a bounded operator. Hence its inverse $\lambda I_{\tilde X} -A_{{\bf T},op} : {\mathcal D}_{{\bf T},op}\to \tilde X$ is closed, and $A_{{\bf T},op}$ is also closed. Now assume that $(X_*,X)$ is also a dual pair, which means that the $\sigma(X_*,X)$-closed convex hull of every $\sigma(X_*,X)$-compact subset of $X_*$ is $\sigma(X_*,X)$-compact. Set $u_{\lambda}(t)=e^{-\lambda t}$ for $\lambda \in \C, t \ge 0.$ Then $R_{{\bf T},op}(\lambda)=\phi_{\bf T}(u_\lambda)$ for $Re(\lambda)> \rho_{\bf T},$ and so it follows from proposition 3.2(iii) that $R_{{\bf T},op}(\lambda)$ is weakly continuous with respect to $X_*,$ which shows as above that $A_{{\bf T},op}$ is weakly closed.

{\color{red} Check that the infinitesimal generator is not closed if lim sup$_{t\to 0^+}\Vert T(t)\Vert =+\infty.$}

\section{Quasimultipliers on the Arveson ideal}

 
  If ${\mathcal U}$ be a commutative Banach algebra, set $\Omega({\mathcal U}):=\{ u \in {\mathcal U} \ | \ \left [ u{\mathcal U}\right ]^-={\mathcal U}\},$ and set ${\mathcal U}^{\perp} =\{ u \in {\mathcal U} \ | \ u{\mathcal U}=\{0 \}\}.$ We recall the definition of quasimultipliers, a notion introduced by the author in \cite{e2}.
  \begin{defn} \mbox{\cite{e2}} Assume that $\Omega({\mathcal U})\neq \emptyset$ and that ${\mathcal U}^{\perp}=\{0\}.$ A quasimultiplier on ${\mathcal U}$ is a pair $(S_{u/v},{\mathcal D}_{S_{u/v}}),$ where $u \in {\mathcal U}, v\in \Omega({\mathcal U}),$ where ${\mathcal D}_{S_{u/v}}$ is the ideal of ${\mathcal U}$ consisting of those $x \in {\mathcal U}$ such that $ux \in v{\mathcal U},$ and where $S_{u/v}x$ is the unique $y \in {\mathcal U}$ such that $ux=vy$ for $x \in {\mathcal D}_{S_{u/v}}.$  The set of quasimultipliers on $\mcu$ will be denoted by ${\qm}(\mcu)$, and a set $B\subset \qm(\mcu)$ will be said to be pseudobounded if there exists $x \in \left [\cap_{S \in B}{\mathcal D}_S\right ] \cap \Omega(\mcu)$ such that sup$_{S\in B}\Vert Sx\Vert <+\infty.$
  \end{defn}
  
 The quasimultipliers on $\mcu$ form an algebra, which is isomorphic to the algebra of fractions $\mcu/\Omega(\mcu)$.  Notice that since $\Omega(\mcu)$ is stable under products, the product of two pseudobounded sets is pseudobounded, which gives
  $\qm(\mcu)$ a structure of "alg\`ebre \`a born\'es".
  
  The algebra $\qm(\mcu)$ is in some sense too large, since $u$ in invertible in $\qm(\mcu)$ for every $u\in \Omega(\mcu),$ so it is natural to consider the following class.
  
  \begin{defn} \mbox{\cite{e2}} A quasimultiplier $R \in \qm(\mcu)$ is said to be regular if there exists $\lambda >0$ such that the family $(\lambda^nR^n)_{n\ge 1}$ is pseudobounded.  The set of regular quasimultipliers on $\mcu$ will be denoted by ${\qm}_r(\mcu)$, and a set $B\subset \qm_r(\mcu)$ will be said to be multiplicatively pseudobounded if it is contained in some set of the form $\lambda V$, where $\lambda >0$ and where $V$ is a pseudobounded subset of $\qm_r(\mcu)$ which is stable under products. 
  
  \end{defn}
  
  Recall that a multiplier on ${\mathcal U}$ is a bounded linear map $S: {\mathcal U} \to {\mathcal U}$ such that $(Su)v=S(uv)$ for $(u,v)\in {\mathcal U}.$ The set ${\mathcal M}({\mathcal U})$ of multipliers on ${\mathcal M}({\mathcal U})$ is a closed subalgebra of ${\mathcal B}({\mathcal U}),$ and ${\mathcal U}$ can be identified to an ideal of  ${\mathcal M}({\mathcal U})$ in an obvious way. We also have the obvious identification ${\mathcal M}({\mathcal U}):=\{ S=S_{u/v}\in \qm({\mathcal U}) \ | {\mathcal D}_{S_{u/v}})={\mathcal U}\}\subset \qm_r({\mathcal U}).$
  
   Equipped with the family of all its multiplicatively pseudobounded subsets, the set $\qm_r(\mcu)$ forms an algebra which is a pseudo-Banach algebra in the sense of Allan, Dales and McClure \cite{ad}. In particular every maximal ideal of  $\qm_r(\mcu)$ is the kernel of a character of  $\qm_r(\mcu),$ and $\qm_r(\mcu)$ is an inductive limit of commutative Banach algebras. In fact regular quasimultipliers on $\mcu$ can be turned into multipliers in the usual sense  by modifying the algebra ${\mathcal U},$ as will be seen later.
  
   Precise estimates on inverse Laplace transforms, obtained by the author in section 6 of \cite{e1} in collaboration with Paul Koosis, combined with Nyman's characterization of dense ideals of $L^1(\R^+),$ show that if $\omega$ is a positive, continuous, nonincreasing submultiplicative weight on $(0,+\infty)$ the convolution algebra $L^1_{\omega}:=L^1_{\omega}(\R^+)$ contains an infinitely differentiable convolution semigroup $(a(t))_{t>0}$ such that $a(t)*L^1_{\omega}$ is dense in $L^1_{\omega}$ for every $t>0$, \cite{e1}, theorem 6.8. This result uses the fact that if $f$ is a continuous function on $(0,+\infty)$ such that $\lim \sup_{t\to \infty}f(t)<+\infty,$ then there exists a continuously differentiable function $h$ on $(0,+\infty)$ such that $h(t)\ge f(t)$ for $t>0,$ 
   and the rather sophisticated contruction, applied to the function $f:t \to e^{-\lambda t}\omega(t)$ for some suitable  $\lambda\in \R,$  consists in finding a convolution semigroup $(a(t))_{t>0}$ of functions in $L^1(\R^+)$ such that $\lim \sup _{t\to 0^+}h(t)\Vert a(t)\Vert <+\infty$ and such that  the convolution ideal generated by $a(t)$ is dense in $L^1(\R^+)$ for every $t>0$ (which implies that the Laplace tranform of $a(t)$ is an outer function of the right-hand half-plane). The fact that $a(t)*L^1_{\omega}$ is dense in $L^1_{\omega}$ follows from Cohen's factorization theorem for modules over commutative Banach algebras with bounded approximate identities, which implies that every $f\in L^1_{\omega}$ can we written under the form $f=g*h,$ where $g \in L^1(\R^+)$ and $h\in L^1_{\omega}.$

   Now assume that $f$ is positive and lower semicontinuous on $(0, +\infty),$ so that $f$ is bounded on $[a,b]$ for $0<a<b<+\infty,$ and that $\lim \sup_{t\to +\infty}f(t)<+\infty.$ Set $f_1(t)=sup_{s\ge t}f(s),$ $f_2(t)={2\over t}\int_{t/2}^tf_1(s)ds.$ Then $f_1$ is nonincreasing, $f_1(s/2) \ge f_2(s)\ge f_1(s)\ge f(s),$ and applying the previous result to $f_2$ we see that there exists a continuously differentiable function $h$ on $(0,+\infty)$ such that $h(t)\ge f(t)$ for $t>0.$ This shows that the construction used in the proof of theorem 6.8 of \cite{e1} extends to lower-semicontinuous weights on $(0,+\infty).$ Set $u_{\lambda}(t)=e^{\lambda t}$ for $t>0.$ Since the map $f\to fu_{-\lambda}$ is a norm-preserving homomorphism from the convolution algebra $L^1_{\omega}(\R^+)$ onto the convolution algebra $L^1_{u_{\lambda}\omega}(\R^+)$ for $\lambda \in \R,$ we obtain the following result.
   
   \begin{thm} Let $\omega$ be a positive, lowersemicontinuous, submultiplicative weight on $(0,+\infty)$, and assume that $\tilde \omega: t \to e^{\lambda t}\omega(t)$ is nonincreasing  for some $\lambda \in \R.$ Then the convolution algebra $L^1_{\omega}:=L^1_{\omega}(\R^+)$ contains an infinitely differentiable semigroup $(a(t))_{t>0}$ such that $a(t)*L^1_{\omega}$ is dense in $L^1_{\omega}$ for every $t>0$.
   
   \end{thm}
   
   We obtain the following  result, which shows in particular that the ideal ${\mathcal I}_{\bf T}$ possesses dense principal ideals, and that the family $(e^{-\lambda t}T(t))_{t>0}$ is pseudobounded in $\qm(\iit)$ for $\lambda >\lim_{t\to +\infty}\Vert T(t)\Vert^{1\over t}.$ We set again  $\rho_{\bf T}:=\lim_{t\to +\infty}\Vert T(t)\Vert^{1\over t}.$
   
   \begin{prop} Let $\lambda > log\left (\rho_{\bf T}\right )$,  and set $\omega_{\lambda}(t)=e^{\lambda t }\sup_{s\ge t}e^{-\lambda s}\Vert T(s)\Vert$ for $\lambda>0.$ Let $g \in L^1_{\omega_{\lambda}}.$  Then $\Vert \phi_{\bf T}(g)T(t)\Vert \le e^{\lambda t}\Vert g \Vert_{L^1_{\omega_{\lambda}}}$ for every $g \in L^1_{\omega_{\lambda}}$ and every $t>0,$ and $\phi_{\bf T}(g) \subset \Omega(\iit)$ for every $g \in \Omega \left (L^1_{\omega_{\lambda}}\right ).$
   \end{prop}
   
   Proof: Notice that $\Vert T(t) \Vert \le \omega_{\lambda}(t)$ for $t>0,$ so that  $L^1_{\omega_{\lambda}}\subset L^1_{\omega_{\bf T}}$. Also $\omega_{\lambda}(s) \ge e{\lambda s}\sup_{t\ge 0}e^{-\lambda (s+t)}\Vert T(s+t)\Vert,$ and so $\Vert T(s+t)\Vert \le e^{\lambda t}\omega_{\lambda}(s)$ for $s\ge 0, t>0.$ We have, for $t>0,$ $g \in L^1_{\omega_{\lambda}},$
   
   $$\Vert \phi_{\bf T}(g)T(t)\Vert =\Vert \phi_{\bf T}(g*\delta_t)\Vert  \le \left \Vert \int_t^{+\infty}g(s-t)T(s)ds \right \Vert \le \int_t^{+\infty}\vert g(s-t)\vert \Vert T(s)\Vert ds$$ $$\le  \int_0^{+\infty}\vert g(s)\vert \Vert T(s+t)\Vert ds \le e^{\lambda t}\Vert g \Vert_{L^1_{\omega_{\lambda}}}.$$
   
   Since $L^1_{\omega_{\bf T}}$ contains $C_c^{\infty}[(0,+\infty)],$ $L^1_{\omega_{\lambda}}$ is dense in $L^1_{\omega_{\bf T}},$ and so $\Omega(L^1_{\omega_{\lambda}}) \subset \Omega(L^1_{\omega_{\bf T}}),$ which shows that $\phi_{\bf T}\left ( \Omega(L^1_{\omega_{\lambda}}) \right )\subset \Omega(\iit).$ $\square$

  \section{Quasimultipliers on the Arveson ideal as weakly densely defined operators}
  
  We now associate to every quasimultiplier on $\iit$ a weakly densely defined operator on $X_1.$ Set again $\tilde X:=\{ x\in X_1 \ | \ \lim_{t\to 0^+} <T(t)x -x,l>=0 \ \forall l \in X_*\}.$
  
  \begin{lem} Let $v \in \Omega(\iit).$ Then $Ker(v) \subset \cap_{t>0}Ker(T(t),$ and $v:\tilde X \to X_1$ is one-to-one.
  \end{lem}
  
  Proof: Let $x \in Ker(v),$ and let $(f_n)_{n\ge 1} \subset {\mathcal C}_c((0,\infty))$ be a Dirac sequence. Since $v\iit$ is dense in $\iit,$ there exists a sequence $(w_n)_{n\ge 1}$ of elements of $\iit$ such that $\lim_{n\to +\infty}\Vert \theta_{\bf T}(f_n)-vw_n\Vert=0.$ It follows then from proposition 2.4 (ii) that we have, for $x \in X, l \in X_*,$ $t>0,$
 
    $$lim_{n\to +\infty}< vT(t)w_nx-T(t)x,l>=0.$$
    
    Since $vT(t)w_nx= T(t)w_nvx=0,$ we have $x \in \cap_{t>0}Ker(T(t)).$ Since $(\cap_{t>0}Ker(T(t))\cap \tilde X=\{0\},$ this shows that $v$ is one-to-one on $\tilde X.$ $\square.$
  
  \begin{prop} Let $S =S_{u/v}\in \qm(\iit).$ Set ${\mathcal D}_{\tilde S}:= \{ x \in  X \ | \ ux\in v(\tilde X)\},$ and denote by $\tilde S(x)$ the unique $y \in \tilde X$ satisfying $ux=vy.$
  Then ${\mathcal D}_{\tilde S}$ is weakly dense in $X_1.$ If lim sup$_{t\to 0^+}\Vert T(t)\Vert <+\infty,$ then $\tilde S$ is a closed partially defined operator on $X_1,$ and 
  $\tilde S$ is weakly closed if, further,  the $\sigma(X_*,X)$-closed convex hull of every $\sigma(X_*,X)$-compact subset of $X_*$ is $\sigma(X_*,X)$-compact.
\end{prop}
  
  Proof:  It follows from the lemma that $\tilde S$ is well-defined on ${\mathcal D}_{\tilde S}.$ Since $v\iit$ is dense in $\iit,$ we see as above there exists a sequence $(w_n)_{n\ge 1}$ of elements of $\iit$ such that  we have, for $x \in X, l \in X_*,$ $t>0,$
  
  $$lim_{n\to +\infty}< vT(t)w_nx-T(t)x,l>=0.$$
  
Since $uvT(t)w_n=vuT(t)w_n\in v(\tilde X),$ we have $T(t)w_nx\in {\mathcal D}_{\tilde S}$ for $t>0, n\ge 1, x \in X,$ and we see that the closure of ${\mathcal D}_{\tilde S}$ in $X_1$ with respect to the weak topology $\sigma(X,X_*)$ contains $T(t)(X)$ for every $t>0,$ which shows that  ${\mathcal D}_{\tilde S}$ is weakly dense in $X_1.$

If lim sup$_{t\to 0^+}\Vert T(t)\Vert <+\infty,$ then $\tilde X=X_1,$ and ${\mathcal D}_{\tilde S}=\{ x \in X_1 \ | \ ux \in v(X_1)\}.$ Let $(x,y) \in X_1\times X_1,$ and assume that there exists
a sequence $(x_n)_{n\ge 1}$ of elements of ${\mathcal D}_{\tilde S}$ such that lim$_{n\to +\infty}\Vert x-x_n\Vert =$lim$_{n\to +\infty}\Vert y-\tilde S x_n\Vert=0.$ Then $ux=$lim$_{n\to +\infty}ux_n=$lim$_{n\to +\infty}vy_n=vy.$ Hence $x \in {\mathcal D}_{\tilde S}$ and $y =\tilde S x,$ which shows that $\tilde S$ is a closed operator. Now if the $\sigma(X_*,X)$-closed convex hull of every $\sigma(X_*,X)$-compact subset of $X_*$ is $\sigma(X_*,X)$-compact, then $u$ and $v$ are weakly continuous, and a similar argument using generalized sequences shows that $\tilde S$ is weakly closed. $\square$

Assume that $\cap_{t>0}Ker(T(t)=\{0\}.$ Then $v$ is one-to-one on $X_1$ for every $v\in \Omega(\iit).$ If $S=S_{u/v} \in \qm(\iit),$ then we can set

\begin{equation} {\mathcal D}_{\overline{\tilde S}}:=\{ x \in X_1 \ | ux \in v(X_1)\},\end{equation},

and define $\overline{\tilde S}: {\mathcal D}_{\overline{\tilde S}}$ by the formula

\begin{equation} v(\overline{\tilde S}x)=ux \ \ \ (x \in {\mathcal D}_{\overline{\tilde S}}).\end{equation}

Then $\overline{\tilde S}$ is a closed extension of $\tilde S,$ and $\overline{\tilde S}$ is also weakly closed if lim sup$_{t\to 0^+}\Vert T(t) \Vert <+\infty.$

To deal with the general case observe that if $x \in X_1,$ and if $T(s)x\in \cap_{t>0}T(t)(X)$ for every $s>0,$ then $T(t)(x)=T(t/2)(T(t/2)x)=0$ for every $t>0.$ So if we set  
$X_0:=X_1/\cap_{t>0}Ker(T(t)),$ and if we set $\tilde T(t)x+\cap_{t>0}Ker(T(t))=T(t)x +\cap_{t>0}Ker(T(t))$ for $x \in X_1, t>0,$ then $\cap_{t>0}Ker(\tilde T(t))=\{0\}.$ Now if $u \in \at,$ then $u(\cap_{t>0}Ker(T(t))\subset  (\cap_{t>0}Ker(T(t)) ,$ and we can define $\tilde u: X_0\to X_0$ by the formula

$$\tilde u(x+(\cap_{t>0}Ker(T(t)))=ux +(\cap_{t>0}Ker(T(t)) ) \ \ (x \in X_1),$$

and the same argument as in the proof of lemma 6.1 shows that if $v\in \Omega(\iit),$ and if $vx \in \cap_{t>0}Ker(T(t))$ for some $x \in X_1,$ then $x \in \cap_{t>0}Ker(T(t)),$ so that $\tilde v$ is one-to-one on $X_0.$

This suggests the following definition.

\begin{defn} Let $S =S_{u/v}\in \qm(\iit).$ Set ${\mathcal D}_{\overline{\tilde S}}:= \{ y \in X_0 /  \ | \ {\tilde u}y\in \tilde v(X_0)\}.$ The operator  $\overline{\tilde S}: {\mathcal D}_{\overline{\tilde S}} \to X_0$ is defined by the formula

$$\tilde v \left (\overline{\tilde S}y\right )= \tilde u y, \ \ (y \in {\mathcal D}_{\overline{\tilde S}}).$$

\end{defn}

We see that $\overline{\tilde S}$ is a closed operator on $X_0,$ which may be considered as an extension of $\tilde S$ since $\left (\cap_{t>0}Ker(T(t))\right )\cap \tilde X=\{0\}.$

  \section{The normalized Arveson ideal}
  
  It follows from theorem 5.3 and proposition 5.4 that   the family $(e^{-\lambda t}T(t))_{t>0}$ is pseudobounded in $\qm(\iit)$ for $\lambda >\lim_{t\to +\infty}\Vert T(t)\Vert^{1\over t}.$We now apply to the Arveson ideal $\iit$ and to the family $(T(t))_{t>0}$ a slight modification of construction performed by the author in \cite{e2} to embed a pseudobounded subset $U$ stable under products of the quasimultiplier algebra $\qm({\mathcal A})$ of a Banach algebra ${\mathcal A}$ such that ${\mathcal A}^{\perp}=\{0\}$ and $\Omega({\mathcal A})\neq \emptyset$ into a bounded set of the multiplier algebra of some Banach algebra ${\mathcal B}$ which is "similar" to ${\mathcal A}.$ In order to do this set ${\mathcal L}:=\{ u \in {\mathcal A} \ | \ \Vert u\Vert_{\mathcal L}:=sup_{b \in U}\Vert bu\Vert <+\infty\},$ choose $b\in  {\mathcal L}\cap\Omega({\mathcal A}),$ and denote by  ${\mathcal U}$ the closure of $b{\mathcal L}$ in $(\mathcal L, \Vert .\Vert_{\mathcal L});$ then $b^2 \in \Omega({\mathcal U}).$ The algebra ${\mathcal B}$ introduced in \cite{e2} is the closure of ${\mathcal A}$ in ${\mathcal M}({\mathcal L}),$ the map $T{u/v} \to T_{b^2u/b^2v}$ defines a pseudobounded isomorphism from $\qm(\mathcal B)$ onto $\qm(\mathcal A)$ and $\qm(\mathcal U),$ and the given pseudobounded set $U$ is contained in the unit ball of ${\mathcal M}(\mathcal B).$ This construction was inspired by Feller's discussion \cite{f1} of the generator of strongly continuous semigroups which are not bounded near the origin, and we slightly modify it below by replacing the multiplier algebra ${\mathcal M}({\mathcal L})$ by the multiplier algebra ${\mathcal M}({\mathcal U})$ in order to obtain a norm-preserving  isomorphism between the multiplier algebras ${\mathcal M}({\mathcal U})$ and  ${\mathcal M}({\mathcal B})$ in the case where $U=\{e^{-\lambda t}T(t)\}$ for some sufficiently large real number $\lambda.$ We will use the convention $T(0)u=u$ for $u\in \iit.$
  
  \begin{thm} Set ${\mathcal L}_{\bf T}:=\{ u \in \iit \ | \ \lim \sup_{t\to 0^+}\Vert T(t)u\Vert <+\infty\}\supset \cup_{t>0}T(t)\iit,$ choose $\lambda > log(\rho_{\bf T}),$ set $\Vert u \Vert_{\lambda}:=sup_{s\ge 0}e^{-\lambda s}\Vert T(s)u\Vert$ for $u \in {\mathcal L}_{\bf T},$ denote by  ${\mathcal U}_{\bf T}$ the closure of  $\cup_{t>0}T(t){\mathcal I}_{\bf T}$ in $({\mathcal L}_{\bf T}, \Vert . \Vert_{\lambda}),$ and set  $\Vert R \Vert_{\lambda,op} := \sup\{\Vert Ru\Vert _{\lambda} \ | \ u \in {\mathcal U}_{\bf T}, \Vert u \Vert_{\lambda}\le 1\}=\Vert R \Vert_{{\mathcal M}({\mathcal U}_{\bf T})}$ for $R\in {\mathcal M}({\mathcal U}_{\bf T}).$ 
  
  Then $({\mathcal L}_{\bf T}, \Vert .\Vert_{\lambda})$ is a Banach algebra, the norm topology on ${\mathcal L}_{\bf T}$ does not depend on the choice of $\lambda,$ ${\mathcal U}_{\bf T}=\{ u \in {\mathcal L}_{\bf T} \ | \ \lim_{t\to 0^+} \Vert T(t)u -u\Vert_{\lambda}=0\},$ ${\mathcal U}_{\bf T}$ is a dense ideal of $\iit,$  and the following properties hold
  
  (i) $\Vert T(t)u\Vert_{\lambda} \le \left (\sup_{t\ge 0}e^{-\lambda s}\Vert T(t+s)\Vert\right )\Vert u \Vert$ for $u\in \iit, t>0.$
  
  (ii) $\Vert T(t)\Vert_{\lambda,op}\le e^{\lambda t}$ for every $t>0.$
  
 (iii) $R{\mathcal L}_{\bf T}\subset {\mathcal L}_{\bf T}$ and $R({\mathcal U}_{\bf T}) \subset {\mathcal U}_{\bf T}$ for every $R \in {{\mathcal M}(\iit)},$ and $\Vert R\Vert_{\lambda,op} \le \Vert R\Vert_{{\mathcal M}(\iit)}.$  In particular $\Vert u \Vert_{\lambda, op}\le \Vert u \Vert$ for every $u \in \at.$
 
 (iv) $ab\in {\mathcal U}_{\bf T}$ for $a\in {\mathcal L}_{\bf T},$ $b \in {\mathcal L}_{\bf T}.$

 (v) $\Omega(\iit)\cap{\mathcal L}_{\bf T} \neq \emptyset,$ $a{\mathcal L}_{\bf T}$ is dense in ${\mathcal U}_{\bf T}$ for every $a \in \Omega( \iit)\cap{\mathcal L}_{\bf T},$ and $ab\in\Omega({\mathcal U}_{\bf T})$ for $a\in \Omega( \iit)\cap{\mathcal L}_{\bf T}, b\in \Omega(\iit)\cap{\mathcal L}_{\bf T}.$
 
 \smallskip
 
 Now denote by $\jjt$ the closure of $\iit,$ identified to a subset of ${\mathcal M}(\dt),$ in $( {\mathcal M}(\dt), \Vert .\Vert_{\lambda, op}).$  Then $(\dt)$ is a dense ideal of $(\jjt, \Vert .\Vert_{\lambda, op})$ and the following properties hold 
  
  (vi) $(\phi_{\bf T}(f_n*\delta_{\epsilon_n}))_{n\ge 1}$ is a bounded sequential approximate identity for $\jjt$ for evey Dirac sequence $(f_n)_{n\ge 1}$ and every sequence $(\epsilon_n)_{n\ge 1}$ such that $\lim_{n\to +\infty}\epsilon_n=0,$ and $\lim\sup_{n\to +\infty}\Vert \phi_{\bf T}(f_n*\delta_{\epsilon_n})\Vert_{\lambda,op}=1.$ 
  
  (vii) $Ru\in \dt$ for every $R \in {\mathcal M}(\jjt)$ and every $u \in \dt,$ and the map $R \to R_{\vert_{\dt}}$ is a norm-preserving isomorphism from ${\mathcal M}(\jjt)$ onto ${\mathcal M}(\dt).$
  
(viii) If $a \in \Omega(\dt)$ then  the map $T_{u/v}\to T_{au/av}$ defines a pseudbounded isomorphism from $\qm(\jjt)$ onto $\qm(\iit)$ and from $\qm(\jjt)$ onto $\qm(\dt).$

\smallskip
The Banach algebra $\jjt$, which is a closed ideal of its multiplier agebra ${\mathcal M}(\jjt),$ will be called the normalized Arveson ideal of the semigroup ${\bf T}.$
  \end{thm}
  
  Proof: Let $\lambda  > log(\rho_{\bf T}),$ and fix $a>0.$ Set $k(a)=sup_{s\ge a}e^{-\lambda s}\Vert T(s)\Vert.$ Property (i) follows from the definition of the norm $\Vert .\Vert_{\lambda},$ and we have
  
  $$\Vert u \Vert _{\lambda}\le max (\sup_{0\le s \le a}\Vert e^{-\lambda s}T(s)u\Vert, k(a)\Vert u\Vert)\le max(e^{\vert \lambda \vert a}\sup_{0\le s \le a}\Vert T(s)u\Vert, k(a)\Vert u\Vert)$$ $$\le max(e^{\vert \lambda \vert a},k(a))\sup_{0\le s \le a}\Vert T(s)u\Vert,$$
  
  $$\sup_{0\le s \le a}\Vert T(s)u\Vert \le e^{\vert \lambda \vert a}\sup_{0\le s \le a}e^{-\lambda s}\Vert T(s)u\Vert \le  e^{\vert \lambda \vert a}\Vert u \Vert _{\lambda},$$
  
  which shows that all the norms $\Vert . \Vert_{\lambda}$ are equivalent to the norm $u \to \sup_{0\le s \le a}\Vert T(s)u\Vert.$

  The construction of ${\mathcal L}_{\bf T}$  is the first part of the construction of theorem 7.11 of \cite{e2} applied to the pseudobounded family stable under products $\{e^{-\lambda t}T(t)\}_{t>0}.$ It follows then from theorem 7.11 of \cite{e2} that $({\mathcal L}_{\bf T}, \Vert .\Vert_{\lambda})$ is a Banach algebra, and (ii) follows from the fact that $\Vert T(t)u\Vert_{\lambda}\le e^{\lambda t}\Vert u \Vert_{\lambda}$ for $u \in {\mathcal L}_{\bf T}, t>0.$ It also follows from theorem 7.11 of \cite{e2} that ${\mathcal L}_{\bf T}$ is an ideal of ${\mathcal M}(\iit)\supset \at $ and that we have the obvious inequality $\Vert R \Vert_{op,\lambda} \le \Vert R \Vert_{{\mathcal M}(\iit)}$ for every $R \in {\mathcal M}(\iit),$ which implies that ${\mathcal L}_{\bf T}$ is an ideal of $\at$ and that $\Vert a \Vert_{op, \lambda}\le \Vert a \Vert$ for every $a \in \at.$ This implies that $R(\dt) \subset \dt$ for every $R \in {\mathcal M}(\iit),$ which completes the proof of (iii).
  
  Set $q(u)=\sup_{0\le s\le 1}\Vert T(s)u\Vert$ for $u \in \lt,$ so that the norm $q$ is equivalent to the norm $\Vert .\Vert_{\lambda}$ on $\lt.$
  
  Let $t>0$ and let $u\in \iit.$ Set $k:=sup_{t/2 \le s \le t/2+1}\Vert T(s)\Vert.$ Then $q(T(t+r)u-T(t)(u))\le k \Vert T(t/2+r)u-T(t/2)u\Vert$ and so $\lim \sup_{r\to 0^+}q(T(r)T(t)u-T(t)u)=0.$
  
  It follows then from (ii) that $\tilde \iit\subset \{ u \in \lt \ | \ \lim_{t\to 0^+} \Vert T(t)u -u\Vert_{\lambda}=0\},$ and in fact $\dt= \{ u \in \lt \ | \ \lim_{t\to 0^+} \Vert T(t)u -u\Vert_{\lambda}=0\},$ since the other inclusion is trivial.
  
  Let $a \in \lt,$ $b \in \lt.$ There exists a sequence $(u_n)_{n\ge 1}$ of elements of $\cup_{t>0}T(t)\iit$ such that $\lim_{n\to +\infty}\Vert b -u_n\Vert=0,$ and we have
  
  $$\lim \sup_{n\to +\infty}\Vert  ab-au_n\Vert _{\lambda}\le \Vert a \Vert_{\lambda} \lim \sup_{n\to +\infty}\Vert b-u_n\Vert =0.$$ Hence $ab \in \dt,$ which proves (iv).
  
  Now let $a \in \Omega(\iit)\cap \lt.$ Then $a\lt \subset \dt.$ Let $u\in \iit,$ and let $(u_n)_{n\ge 1}$ be a sequence of elements of $\iit$ such that $\lim_{n\to +\infty}\Vert u-au_n\Vert=0.$ It follows from (i) that $\lim_{n\to +\infty} \Vert T(t)u -T(t)au_n\Vert_{\lambda}=0,$ and so the closure of $a\lt$ in $(\dt, \Vert .\Vert)$ contains $\cup_{t>0}T(t)\iit,$ hence coincides with $\dt.$ An easy verification given in the proof of theorem 7.11 of \cite{e2} shows then that $ab \in \Omega(\dt)$ for $a \in \Omega(\iit)\cap \lt, b \in \Omega(\iit)\cap \lt.$ Notice that the ideals $\lt$ and $\dt$ play with respect to $\iit$ and the family $(e^{-\lambda t}T(t)$ the role played by the ideals ${\mathcal L}$ and ${\mathcal U}$ with respect to the Banach algebra $\mathcal A$ and the pseudobounded family $U$ stable under products in the construction of theorem 7.11 of \cite{e2}.
  
    Now identify $\iit$ to a subset of ${\mathcal M}(\dt),$ and denote by $\jjt$ the closure of $\iit$ in $( {\mathcal M}(\dt), \Vert .\Vert_{\lambda, op}).$ By construction, $\dt$ is an ideal of $\jjt$ which is dense in $\jjt$ since it contains $\cup_{t>0}T(t)\iit.$ Then  the Banach algebras $\dt, \iit$ and $\jjt$ are similar in the sense of \cite{e2}, definition 7.4 and the natural injection $\iit \to \jjt$ in a $s$-homomorphism in the sense of \cite{e2}, definition 7.9, and it follows from theorem 7.11 of \cite{e2} that (viii) holds.
   
   Let $(f_n)_{n\ge 1}$ be a Dirac sequence and let $(\epsilon_n)_{n\ge 1}$ be a sequence of positive real numbers such that $\lim_{n\to +\infty}\epsilon_n=0.$ There exists a sequence $(\delta_n)_{n\ge 1}$ of positive real numbers such that $\lim_{n\to +\infty}\delta_n=0$ and such that $f_n(t)=0$ a.e. for $t >\delta_n.$ Set $e_n=\phi_{\bf T}(f_n*\delta_{\epsilon_n}).$ We have
   
   $$\Vert e_n\Vert_{\lambda,op}\le \int_0^{\delta_n}f_n(t)\Vert T(t+\epsilon_n)\Vert dt,$$
   
   and it follows from (ii) that $\lim \sup _{n\to +\infty}\Vert e_n \Vert_{\lambda,op}\le 1.$ Since $\lim_{n\to +\infty}\Vert e_nu-u\Vert=0$ for $u\in \cup_{t>0}T(t)\iit,$ we have a fortiori $\lim_{n\to +\infty}\Vert e_nu-u\Vert_{\lambda} =0$ for $u\in \cup_{t>0}T(t)\iit,$ and a standard density argument shows that $\lim_{n\to +\infty}\Vert e_nu-u \Vert_{\lambda =0}$ for $u\in \jjt.$ Hence $(e_n)_{n\ge 1}$ is a sequential approximate identity for $\jjt$ and satisfies $lim_{n\to +\infty}\Vert e_n\Vert_{op, \lambda}=1.$ This implies that $\Vert u \Vert_{\mathcal M(\jjt)}=\Vert u \Vert_{\lambda,op}$ for $u\in \jjt.$ 
   
   Let $u \in \jjt,$ let $S\in {\mathcal M}(\dt).$ Since $\dt$ is dense in $(\iit, \Vert.\Vert)$, $\dt$ is dense in $\iit$ with respect to $\Vert.\Vert_{\lambda,op},$ and there exists a sequence $(u_n)_{n\ge 1}$ of elements of $\dt$ such that $\lim_{n\to +\infty}\Vert u -u_n\Vert_{\lambda,op}=0.$ Then $\lim _{n\to +\infty}\Vert Su -Su_n\Vert_{\lambda,op}=0,$ which shows that $Su\in \jjt.$ Hence $S:\dt \to \dt$ has an extension $\tilde S$ to $\jjt,$ and $\Vert \tilde S\Vert_{{\mathcal M}(\jjt)}=\Vert S_{\lambda,op}\Vert.$
   
   Writing $T(t)u=T(t/2)T(t/2)u$ for $t>0, u\in \iit,$ we deduce from proposition 2.4 and $(ii)$ that $\lim_{n\to +\infty}\Vert T(t)u-T(t)ue_n\Vert_{\lambda}=0.$ A standard density argument shows then that $\lim_{n\to +\infty}\Vert ue_n-u\Vert_{\lambda}=0$ for every $u \in \dt.$ It follows then from Cohen's factorization theorem \cite{co} for Banach modules over commutative Banach algebras with bounded approximate identities that for every $u \in \dt$  there exists $v \in \jjt$ and $w \in \dt$ such that $u=vw, \Vert v\Vert_{\lambda,op}=1$ and $\Vert w\Vert _{\lambda}\le \Vert u \Vert _{\lambda +\epsilon}.$ Hence $Ru=(Rv)w\in \dt$ for $R\in {\mathcal M}(\jjt).$ Also $\Vert Ru\Vert_{\lambda} \le \Vert Rv\Vert_{\lambda, op}\Vert \Vert w\Vert_{\lambda} \le \Vert R\Vert_{{\mathcal M}(\jjt)}(\Vert u\Vert_{\lambda}+\epsilon)$ for every $\epsilon >0,$ $\Vert Ru\Vert_{\lambda}\le \Vert R \Vert_{{\mathcal M}(\dt)}\Vert u\Vert_{\lambda},$  which shows that $R_{\l_{\dt}}: \dt \to \dt$ is bounded. So the map $R \to R_{|_{\dt}}$ is a norm-preserving isomorphism from ${\mathcal M}(\jjt)$ onto ${\mathcal M}(\dt).$ $\square$
   
   \section{The generator as a quasimultiplier on the Arveson ideal}
   
   When $S=S_{u/v}$, with $u \in {\mathcal A},$ $v\in \Omega({\mathcal A})$ is a quasimultiplier on a commutative Banach algebra ${\mathcal A}$ such that $\Omega({\mathcal A})\neq \emptyset$ and $a{\mathcal A}\neq \{0\}$ for $a \in {\mathcal A}\setminus\{0\},$ we will write $S=u/v$ if there is no risk of confusion.
    
   \begin{defn}:  The infinitesimal generator $A_{\bf T, \iit}$ of a Arveson weakly continuous semigroup ${\bf T}$ is the quasimultiplier on $\iit$ defined by the formula
 
$$A_{\bf T, \iit}= {-{\phi_{\bf T}(f_0')/\phi_{\bf T}(f_0)}},$$

where $f_0 \in {\mathcal C}^1([0,+\infty))\cap \Omega \left (L^1_{\omega_{\bf T}}\right )$ satisfies $f_0=0,$ $f_0' \in L^1_{\omega_{\bf T}}.$

\end{defn}
   
Assume that $f_1$ also satisfies the conditions of the definition, and set $f_2=f_0*f_1.$ Since $\Omega \left (L^1_{\omega_{\bf T}}\right )$ is stable under convolution, $f_2 \in  \Omega \left (L^1_{\omega_{\bf T}}\right ),$ and $f'_2=f'_0*f_1=f_0*f'_1\in L^1_{\omega_{\bf T}}$ is continuous. Also $f_2(0)=0,$ and we have

$$\phi_{\bf T}(f'_2)/\phi_{\bf T}(f_2)= \phi_{\bf T}(f'_0)\phi_{\bf T}(f_1)/\phi_{\bf T}(f_0)\phi_{\bf T}(f_1)= \phi_{\bf T}(f_0')/b_0\phi_{\bf T}(f_0),$$

and similarly $\phi_{\bf T}(f'_2)/\phi_{\bf T}(f_2)=\phi_{\bf T}(f'_1)/\phi_{\bf T}(f_1),$ which shows that the definition of $A_{\bf T, \iit}$ does not depend on the choice of $f_0.$

Let $\lambda > log(\rho_{\bf T}),$ and set again $u_{\lambda}(t)=e^{\lambda t}$ and $\omega_{\lambda}(t)=e^{\lambda t}\sup_{s\ge t}e^{-\lambda s}\Vert T(s)\Vert$ for $t>0.$
Then $L^1_{u_{\lambda}}*L^1_{\omega _{\lambda}}\subset L^1_{\omega _{\lambda}},$ and every Dirac sequence $(e_n)_{n\ge 1}$ is a sequential bounded approximate identity for
 $L^1_{u_\lambda}.$ It follows from Cohen's factorization theorem \cite{co} that every $f\in L^1_{\omega_{\bf \lambda}}$ can be written under the form $f=g*h,$ where $g\in L^1_{u_{\lambda}}$ and $h\in L^1_{\omega_{\lambda}},$ which implies that $\Omega(L^1_{u_{\lambda}})\cap L^1_{\omega_{\lambda}}\subset \Omega(L^1_{\omega_{\lambda}}).$
 
 It is well-known that the subalgebra of  the convolution algebra $L^1(\R^+)$ generated by $u_{\alpha}$ is dense in $L^1(\R^+),$ for $\alpha <0,$ see for example \cite {e1}, so that $u_{\alpha}$ in $\Omega(L^1(\R^+))$ (this follows also from Nyman's theorem \cite{ny}). So the function $u_{\alpha}*u_{\alpha}:t\to te^{\alpha t}$ generates a dense ideal of $L^1(\R^+),$
 and if $log(\rho_{\bf T}) \le \mu \le \lambda$ then the function $v_{\lambda}:t \to t e^{-\lambda t}$ generates a dense ideal of $L^1_{u_{\mu}}.$ Notice that $v_{\lambda} \in {\mathcal C^{\infty}([0,+\infty})),$ that $v'_{\lambda}\in L^1_{u_{\mu}}$ and that $v_{\lambda}(0)=0.$ Using proposition 5.4, we obtain the following result.
 
 \begin{prop}: Let $\mu > log(\rho_{\bf T}),$ let $g_0 \in  {\mathcal C}^1([0,+\infty))\cap \Omega \left (L^1_{u_{\mu}}\right )$ such that $g_0(0)=0, g'_0 \in L^1_{u_{\mu}},$ and let 
 $h_0 \in L^1_{\omega_{\mu}}\cap \Omega \left (L^1_{u_{\mu}}\right ).$ Then $\phi_{\bf T}(g_0*h_0) \in \Omega(\iit),$ and 
 
 $$A_{\bf T, \iit}= -{\phi_{\bf T}(g_0'*h_0)/\phi_{\bf T}(g_0*h_0)}.$$
 
 \end{prop}
 
 For example we can take $g_0=v_{\lambda}$ for some $\lambda > log(\rho_{\bf T}),$ so that $g_0(t)=te^{-\lambda t}$ and $g'_0(t)=(1-\lambda t)e^{-\lambda t}$ for $t \ge 0.$

We can also consider the infinitesimal generator of the semigroup as a quasimultiplier on the normalized Arveson ideal $\jjt.$ Set $\tilde \omega _{\bf T}:=\Vert T(t)\Vert_{{\mathcal M}(\jjt)},$ and for $\mu\in {\mathcal M}_{\tilde \omega_{\bf T}},$ define $\phi_{\bf T} (\mu)\in {\mathcal M}(\jjt)$ by the formula

\begin{equation}\tilde \phi_{\bf T} (\mu)u=\int_0^{+\infty}T(t)ud\mu(t)   \ \ \ (u \in \jjt).\end{equation}

Clearly, $\tilde \phi_{\bf T} (\mu)u=\phi_{\bf T} (\mu)u$ for $u \in \iit, \mu \in \mathcal M_{\omega_{\bf T}},$ and so we can consider $\tilde \phi_{\bf T}:  {\mathcal M}_{\tilde \omega_{\bf T}} \to {\mathcal M}(\jjt)$ as an extension to ${\mathcal M}_{\tilde \omega_{\bf T}}$ of $\phi_{\bf T}: {\mathcal M}_{\omega_{\bf T}}\to {\mathcal M}(\iit).$

   \begin{defn}:  The normalized infinitesimal generator $A_{\bf T, \jjt}$ of a Arveson weakly continuous semigroup ${\bf T}$ is the quasimultiplier on $\jjt$ defined by the formula
 
$$A_{\bf T, \jjt}= -{\tilde \phi_{\bf T}(g_0')/\tilde \phi_{\bf T}(g_0)},$$

where $g_0 \in {\mathcal C}^1([0,+\infty))\cap \Omega \left (L^1_{\tilde \omega_{\bf T}}\right )$ satisfies $g_0=0,$ $g_0' \in L^1_{\tilde \omega_{\bf T}}.$

\end{defn}

  Let $w \in \Omega(\dt).$ The map $S_{u/v}\to S_{wu/wv}$ defines a pesudobounded isomorphism from $\qm(\jjt)$ onto $\qm(\iit).$ Since $\tilde \phi_{\bf T} (g_0)\phi_{\bf T}(h_0)=\phi_{\bf T}(g_0*h_0)$ and $\tilde \phi_{\bf T} (g'_0)\phi_{\bf T}(h_0)=\phi_{\bf T}(g'_0*h_0)$ if $g_0$ and $h_0$ satisfy the conditions of proposition 8.2, the infinitesimal generator 
  $A_{\bf T, \iit}$ is the quasimultiplier on $\iit$ associated to $A_{\bf T, \jjt}$ via this isomorphism.
  
  \begin{lem} Let $\omega$ be a lower semicontinuous submultiplicative weight on $(0,+\infty),$ and let $f \in {\mathcal C}^1([0,+\infty))\cap L^1_{\omega}.$ If 
  $f(0)=0,$ and if $f' \in L^1_{\omega},$ then the Bochner integral $\int_t^{+\infty}(f'*\delta_s)ds$ is well-defined in $L^1_{\omega}$ for $t\ge 0,$ and we have
  
  $$f*\delta_t  -f =-\int_0^{t}(f'*\delta_s)ds, \ \ \mbox{and} \lim_{t\to 0^+}\left \Vert \frac{f*\delta_t-f}{t}+f'\right \Vert _{L^1_{\omega}}=0.$$
  
  \end{lem}
  
  Proof: Set $R(s)g=g*\delta_s$ fo $s>0, g \in L^1_{\omega}.$ Then $R(s)$ is a strongly continuous semigroup of bounded operators on $L^1_{\omega},$ and $\Vert R(s)\Vert \le \omega(s)$ for $s>0,$ and so the integral $\int_0^{+\infty}(f'*\delta_s)ds$ is well-defined in $L^1_{\omega}.$ Since $f'$ is continuous, this integral is also well-defined in $L^1_{\omega}\cap {\mathcal C}([0,L])$ for every $L>0,$ and we have, for $x \ge 0,$
  
  $$\left (\int_0^{t}(f'*\delta_s)ds\right )(x)=<\int_0^{t}(f'*\delta_s)ds,\delta_x>=\int_0^{t}<f'*\delta_s,\delta_x>ds$$ $$=\int_0^{t}(f'*\delta_s)(x)ds.$$
  
  We have $(f'*\delta_s)(x)=f'(x-s)$ for $0\le s \le x,$ and $(f'*\delta_s)(x-s)=0$ for $ s> x.$ Hence $\int_0^{t}(f'*\delta_s)(x)ds= \int_0^{\min(x,t)t}(f'(x-s)ds=-f(x-\min(x,t))+f(x).$ But $f(x-\min(x,t))=f(x-t)$ for $x\ge t,$ and $f(x-\min(x,t))=f(0)=0$ for $0\le x \le t.$ So $f*\delta_t  -f =-\int_0^{t}(f'*\delta_s)ds$ for $t \ge 0,$ and $\lim_{t\to 0^+}\left \Vert \frac{f*\delta_t-f}{t}+f'\right \Vert _{L^1_{\omega}}=0.$ $\square$

  We now give the classical approach to the domain of the generator of a semigroup.
  
  \begin{prop} (i) Let $u \in \iit.$ If $\lim_{t\to 0^+}\Vert \frac{T(t)u-u}{t}-v\Vert =0$ for some $v \in \dt,$ then  $u \in {\mathcal D}_{A_{\bf T, \iit}},$ and $A_{\bf T, \iit}u=v.$
  
  (ii) Let $u \in \jjt.$ Then $u \in {\mathcal D}_{A_{\bf T, \jjt}}$ if, and only if, there exists $v \in \jjt$ such that $\lim_{t\to 0^+}\left \Vert \frac{T(t)u-u}{t}-v\right \Vert _{\jjt}=0,$ and in this situation
  $A_{\bf T, \jjt}u=v.$
  
  \end{prop}
  
  Proof: (i) If $u \in \iit,$ and if $\lim_{t\to 0^+}\Vert \frac{T(t)u-u}{t}-v\Vert =0$ for some $v \in \dt,$ let $f_0 \in {\mathcal C}^1([0,+\infty))\cap \Omega \left (L^1_{\omega_{\bf T}}\right )$ satisfiying $f_0=0,$ $f_0' \in L^1_{\omega_{\bf T}}.$ It follows from the lemma that we have in $\iit$
  
  $$-\phi_{\bf T}(f'_0)u=\left [\lim_{t\to 0^+}\frac {T(t)\phi_{\bf T}(f_0)-\phi_{\bf T}(f_0)}{t}\right ]u= \phi_{\bf T}(f_0)\left [\lim _{t\to 0^+} {T(t)u-u\over t}\right ]= \phi_{\bf T}(f_0)v,$$
  
  and so $u \in  {\mathcal D}_{A_{\bf T, \iit}},$ and $A_{\bf T, \iit}u=v.$
  
  (ii) The same argument shows that if $u\in \jjt,$ and if $\lim_{t\to 0^+}\Vert \frac{T(t)u-u}{t}-v\Vert _{\jjt}=0$ for some $v \in \jjt,$ then $u\in {\mathcal D}_{A_{\bf T, \jjt}},$ and $A_{\bf T, \jjt}u=v.$ Conversely assume that $u\in {\mathcal D}_{A_{\bf T, \jjt}},$ and let $f_0 \in {\mathcal C}^1([0,+\infty))\cap \Omega \left (L^1_{\tilde \omega_{\bf T}}\right ).$ Set $v=A_{{\bf T},\jjt}u.$  We have, for $t\ge 0,$
  
  $$\tilde \phi _{\bf T}(f_0)\int_0^tT(s)vds= \int_0^tT(s)\tilde \phi_{\bf T}(f_0)vds=- \int_0^tT(s)\tilde \phi_{\bf T}(f'_0)uds$$ $$=-\left [ \int_0^tT(s)\tilde \phi_{\bf T}(f'_0)ds\right ]u=\left [T(t) \tilde \phi_{\bf T}(f_0)-\tilde \phi_{\bf T}(f_0)\right ]u=\tilde \phi_{\bf T}(f_0)(T(t)u-u).$$
  
  Since $ \tilde \phi _{\bf T}(f_0)\in \Omega(\jjt),$ this shows that $T(t)u-u=\int_0^tT(s)vds,$ and so $\lim_{t\to 0^+}\left \Vert \frac{T(t)u-u}{t}-v\right \Vert _{\jjt}=0.$ $\square$

\section{The Arveson spectrum}

 We will  denote by $\sia$ the space of characters of $\iit,$ equipped with the usual Gelfand topology. Notice that if $\chi \in \sia$ then there exists a unique character
$\tilde \chi$ on ${\mathcal {QM}}(\iit)$ such that $\tilde \chi _{|_{\iit}}=\chi,$ which is defined by the formula $\tilde \chi(S_{u/v})=\frac{\chi(u)}{\chi(v)}$ for $u\in \iit, v \in \Omega(\iit).$ 

\begin{defn} Assume that $\iit$ is not radical, and let $S \in \qm(\iit).$ The Arveson spectrum $\sigma_{ar}(S)$ is defined by the formula

$$\sigma_{ar}(S)=\{\lambda = \tilde \chi(S) \ : \ \chi \in \sia\}.$$

\end{defn}

If $\mu$ is a measure on $[0,+\infty),$ the Laplace tranform of $\mu$  is defined by the usual formula ${\mathcal L}(\mu)(z)=\int_0^{+\infty}e^{-zt}d\mu(t)$
when $\int_0^{+\infty}e^{-Re(z)t}d\vert \mu \vert(t)<+\infty.$

We have the following easy observation.

\begin{prop}  Let $\mu \in \M.$ Then we have, for $\chi \in \sia,$ 

\begin{equation} \tilde \chi \left ( \int_{0}^{+\infty}T(t)d\mu(t)\right)={\mathcal L}(\mu)(-\tilde \chi(A_{{\bf T},\iit})).\end{equation}

Similarly we have, for $\mu \in {\mathcal M}_{\tilde \omega_{\bf T}}, \chi \in \sia,$

\begin{equation} \tilde \chi \left ( \int_{0}^{+\infty}T(t)d\mu(t)\right)={\mathcal L}(\mu)(-\tilde \chi(A_{{\bf T},\iit}))={\mathcal L}(\mu)(-\tilde \chi(A_{{\bf T},\jjt})) .\end{equation}

In particular $\tilde \chi (T(t))=e^{-\tilde \chi(A_{{\bf T},\iit})t}$ for $t>0.$



\end{prop}

Proof: If $\chi \in \sia,$ then $\tilde \chi_{|_{\at}}$ is a character on $\at,$ the map $t \to \tilde \chi (T(t))$ is continuous on $(0,+\infty)$ and so there exists $\lambda \in \C$ such that $\tilde \chi (T(t))=e^{-\lambda t}$ for $t>0,$ and $\vert e^{-\lambda t}\vert\le \Vert T(t)\Vert,$ which shows that $Re(\lambda)\ge -log(\rho_{\bf T}).$

Let $u\in \Omega(\iit),$ and let $\mu \in \M.$ We have

$$\chi (u)\tilde \chi \left (\int_0^{+\infty}T(t)d\mu(t)\right )= \chi \left (u\int_0^{+\infty}T(t)d\mu(t)\right )
=\chi \left (\int_0^{+\infty}T(t)ud\mu(t)\right )$$
$$=\int_0^{+\infty}\chi(T(t)u)d\mu(t)=\chi(u)\int_0^{+\infty}e^{-\lambda t}d\mu(t)=\chi(u){\mathcal L}(\mu)(\lambda),$$

and so $\tilde \chi(\mu)={\mathcal L}(\mu)(\lambda).$

Let $f_0\in {\mathcal C}^1((0,+\infty))\cap \Omega(\iit)$ such that $f_0(0)=0.$ We have

$$\lambda {\mathcal L}(f_0)(\lambda)={\mathcal L}(f'_0)(\lambda)=\chi (\phi_{\bf T}(f'_0))=-\tilde \chi\left (A_{{\bf T},\iit}\phi_{\bf T}(f_0)\right )$$ $$=-\tilde \chi(A_{{\bf T},\iit})\chi (\phi_{\bf T}(f_0)= -\tilde \chi(A_{{\bf T},\iit}){\mathcal L}(f_0)(\lambda),$$

and so $\lambda =-\tilde \chi(A_{{\bf T},\iit}),$ which proves (14), and  formula (15) follows from a similar argument. In particular $\chi(T(t))={\mathcal L}(-\tilde \chi(A_{{\bf T},\iit}))$ for $t>0.$ $\square$




  The following consequence of  proposition 9.2 pertains to folklore.

\begin{cor} Assume that  $\iit$ is not radical. Then the map $\chi \to \tilde \chi(A_{{\bf T},\iit})$ is a homeomorphism from $\widehat{\iit}$ onto $\sigma_{ar}(A_{{\bf T},\iit}),$ and the set
$\Delta_{t}:=\{ \lambda \in \sigma_{ar}(A_{{\bf T},\iit}) \ | \ Re(\lambda)\le t\}$ is compact for every $t \in \R.$
\end{cor}

Proof: Let $f_0\in {\mathcal C}^1((0,+\infty))\cap \Omega(\iit)$ such that $f_0(0)=0.$ We have $\chi(\phi_{{\bf T}}(f_0))\neq 0$ and $\tilde \chi(A_{{\bf T},\iit})=-\frac{\chi(\phi_{\bf T}(f'_0))}{\chi(\phi_{\bf T}(f_0))}$ for $\chi \in \widehat{\iit},$ 
and so  the map $\chi \to \tilde \chi(A_{\bf T})$ is continuous with respect to the Gelfand topology on $\widehat{\iit}.$

Conversely let $f \in \LL.$ It follows from proposition 5.6 that we have, for $\chi \in \widehat{\iit},$

$$\chi(\phi_{\bf T}(f))={\mathcal L}(f)(-\tilde \chi(A_{{\bf T},\iit})).$$

Since the set $\{ u= \phi_{\bf T}(f) : f \in \LL \}$ is dense in $\iit,$ this shows that the map $\chi \to \tilde {\chi}(A_{{\bf T},\iit})$ is one-to-one on  $\widehat{\iit},$ and that the inverse map $\sigma_{ar}(A_{{\bf T},\iit})\to \widehat{\iit}$ is continuous with respect to the Gelfand topology.

Now let $t \in \R,$ and set $U_t:=\{ \chi \in \widehat{\iit} :  Re(\chi(A_{{\bf T},\iit}))\le t\}.$ Then $\vert \tilde {\chi}(T(1))\vert\ge e^{-t}$ for $\chi \in U_t,$ and so $0$ does not belong to the closure of $U_t$ with respect to the weak$*$ topology on the unit ball of the dual of $\iit.$ Since $\widehat{\iit}\cup \{0\}$ is compact with respect to this topology, $U_t$ is a compact subset of $\widehat{\iit}$, and so the set $\Delta_t$ is compact. $\square$

\section{The resolvent}
We now wish to discuss the Arveson resolvent of the generator of a semigroup ${\bf T}=(T(t))_{t>0}$ which is weakly continuous in the sense of Arveson. We identify $\jjt, {\mathcal M}(\jjt)$ and $\at$ to  subalgebras of the algebra ${\mathcal QM}_r(\iit)$ of regular quasimultipliers on the Arveson ideal $\iit.$ We also identify the algebras ${\mathcal QM}(\iit)$ and ${\mathcal QM}(\jjt),$
and the algebras ${\mathcal QM}_r(\iit)$ and ${\mathcal QM}_r(\jjt)$ by using the pseudobounded isomorphisms described in the previous sections. From now on we will write $A_{\bf T}= 
A_{{\bf T},\iit}=A_{{\bf T},\jjt},$  and we will denote as before by $\mathcal D_{A_{\bf T},\jjt}$ (resp. $\mathcal D_{A_{\bf T},\iit}$, resp. $\mathcal D_{A_{\bf T},\dt})$ the domain of $A_{\bf T}$ considered as a quasimultiplier on $\jjt$ (resp. on $\iit,$ resp. on $\dt$).

\begin{defn} The Arveson resolvent set of  $A_{\bf T}$  is defined by the formula
$$Res_{ar}(A_{\bf T})=\C \setminus \sigma_{ar}(A_{{\bf T}}),$$

with the convention $\sigma_{ar}(A_{{\bf T}})=\emptyset$ if $\iit$ is radical.
\end{defn}

We will denote by $I$ the unit element of ${\mathcal QM}(\iit),$ so that $\jjt \oplus \C I$ is a closed subalgebra of the Banach algebra ${\mathcal M}(\jjt).$ We set as above $u_{\lambda}(t)=e^{\lambda t}$ for $\lambda \in \C, t \ge 0.$ We obtain the usual "resolvent formula".

\begin{prop} The quasimultiplier $\lambda I-A_{\bf T} \in {\mathcal QM}(\iit)$ admits an inverse $(\lambda I-A_{\bf T})^{-1}\in \jjt \subset {\mathcal QM}_r(\iit)$ for $\lambda \in Res_{ar}(A_{\bf T}),$ and the map $\lambda \to (\lambda I -A_{\bf T})^{-1}$ is an holomorphic map from $Res_{ar}(A_{\bf T})$ into $\jjt.$ Moreover we have, for $Re(\lambda) > log(\rho_{\bf T}),$

$$(\lambda I -A_{\bf T})^{-1}=\tilde \phi({\bf T})(u_{-\lambda})=\int_0^{+\infty}e^{-\lambda s}T(s)ds \in \jjt,$$

where the Bochner integral is computed with respect to the strong operator topology on ${\mathcal M}(\jjt),$ and $\Vert \lambda I -A_{\bf T})^{-1}\Vert_{\jjt}\le \int_0^{+\infty}e^{-Re(\lambda)t}\Vert T(t)\Vert dt.$

\end{prop}

Proof: Assume that $Re(\lambda)>log(\rho_{\bf T}),$ let $v \in \jjt,$ and set $a=\tilde \phi_{\bf T}(u_{-\lambda}).$ We have

$$av=\int_0^{+\infty}e^{-\lambda s}T(s)vds, T(t)av-av=\int_0^{+\infty}e^{-\lambda s}T(s+t)vds-\int_0^{+\infty}e^{-\lambda s}T(s)vds$$ $$=e^{\lambda t}\int_t^{+\infty}e^{-\lambda s}T(s)vds -int_0^{+\infty}e^{-\lambda s}T(s)vds= (e^{\lambda t}-1)av -e^{\lambda t}\int_0^te^{-\lambda s}T(s)vds.$$

Since $\lim_{t\to 0^+}\Vert T(t)v-v\Vert_{\jjt}=0,$ we obtain

$$\lim_{t\to 0^+}\left \Vert \frac{T(t)av-av}{t} -\lambda av -v\right \Vert_{\jjt}=0,$$

and so $av \in \mathcal D_{A_{\bf T},\jjt},$ and $A_{\bf T}(av)=\lambda av -v.$ This shows that $a\jjt \subset D_{A_{\bf T},\jjt},$ and that $(\lambda I -A_{\bf T})av=v$ for every
$v \in \jjt.$ We have $\lambda  I -A_{\bf T} =S_{u/v},$ where $u \in \jjt, v\in \Omega(\jjt),$ and we see that $ua=v.$ Hence $u\in \Omega(\iit), \lambda  I -A_{\bf T}$ is invertible in $\qm(\jjt),$ and $(\lambda  I -A_{\bf T})^{-1}=a=\tilde \phi_{\bf T}(u_{-\lambda})=\int_0^{+\infty}e^{-\lambda t}T(t)dt \in \jjt,$ where the Bochner integral is computed with respect to the strong operator topology on ${\mathcal M}(\jjt).$

Set again $a=\tilde \phi_{\bf T}(u_{-\lambda})= (\lambda I-A_{\bf T})^{-1},$ where $\lambda$ is a complex number satisfying $Re(\lambda) > log(\rho_{\bf T}).$ Let $\chi_0$ be the character on $\jjt \oplus \C I$ such that $Ker({\chi_0})=\jjt.$ Every character on $\jjt \oplus \C I$ distinct from $\chi_{0}$ has the form $\tilde \chi_{|_{\jjt \oplus \C I}}$ for some $\chi \in \sia.$ It follows from formula (15) that $\tilde \chi (a)={\mathcal L}(u_{-\lambda})(-\tilde \chi(A_{\bf T}))={1\over \lambda -\tilde \chi(A_{\bf T})}$ for $\chi \in \sia,$ and we obtain, for $\mu \in \C,$

$$spec_{\jjt \oplus \C I}(I+(\mu -\lambda)a)=\{1 \}\cup \left \{ \frac{\mu -\tilde \chi(A_{\bf T})}{\lambda -\tilde \chi(A_{\bf T})}\right \}_{\chi \in \widehat{\iit}}.$$

So $I+(\mu  -\lambda)a$ is invertible in $\jjt \oplus \C e$ for $\mu \in Res_{ar}(A_{\bf T}),$ and the map $\mu \to a(I+(\mu -\lambda)a)^{-1}\in \iit$ is holomorphic on $Res_{ar}(A_{\bf T}).$ We have, for $\mu \in Res_{ar}(A_{\bf T}),$

$$(\mu I -A_{\bf T})a(I+(\mu-\lambda)a)^{-1}=((\mu-\lambda)I +(\lambda -A_{\bf T}))a(I+(\mu-\lambda)a)^{-1}$$ $$=(I +(\mu-\lambda)a)(I+(\mu-\lambda)a)^{-1}=I.$$

Hence $(\mu I-A_{\bf T})$ has an inverse $(\mu I -A_{\bf T})^{-1}\in \jjt$ for $\lambda \in Res_{ar}(A_{\bf T}),$ and the map $\mu \to (\mu I-A_{\bf T})^{-1}=a(I+(\mu-\lambda)a)^{-1}$ is holomorphic on $Res_{ar}(A_{\bf T})$. $\square$

If we consider $A_{\bf T}$ as a quasimultiplier on $\jjt,$ the fact that $(\mu I -A_{\bf T})^{-1}\in \jjt$ is the inverse of $\mu I -A_{\bf T}$ for $\mu \in Res(A_{\bf T})$ means that $(\mu I -A_{\bf T})^{-1}v \in {\mathcal D}_{A_{{\bf T},\jjt}}$ and that $(\mu I-A_{\bf T})\left ((\mu I -A_{\bf T})^{-1}v\right )=v$ for every $v \in \jjt,$ and that if $w \in {\mathcal D}_{A_{{\bf T},\jjt}},$ then 
$(\mu I -A_{\bf T})^{-1}\left ( (\mu I -A_{\bf T})w\right )=w.$ The situation is slightly more complicated if we consider $A_{\bf T}$ as a quasimultiplier on $\iit$ when $\lim \sup_{t \to 0^+}\Vert T(t)\Vert =+\infty.$ In this case the domain ${\mathcal D}_{(\mu I -A_{\bf T})^{-1},\iit}$ of $(\mu I -A_{\bf T})^{-1}\in {\mathcal QM}(\iit)$ is a proper subspace of $\iit$ containing $\dt \supset \cup_{t>0}T(t)\iit,$ and we have $(\mu I -A_{\bf T})^{-1}v \in {\mathcal D}_{A_{{\bf T},\iit}}$ and $(\mu I-A_{\bf T})\left ((\mu I -A_{\bf T})^{-1}v\right )=v$ for every $v \in {\mathcal D}_{(\mu I -A_{\bf T})^{-1},\iit}.$ Also if $w \in {\mathcal D}_{A_{{\bf T},\iit}},$ then $(\mu I -A_{\bf T})w\in {\mathcal D}_{(\mu I -A_{\bf T})^{-1},\iit},$ and we have 
$(\mu I -A_{\bf T})^{-1}\left ( (\mu I -A_{\bf T})w\right )=w.$

In order to interpret $(\lambda I -A_{\bf T})^{-1}$ as a partially defined operator on $\iit$ for $Re(\lambda) >log(\rho_{\bf T}),$ we can use the formula

\begin{equation}(\lambda I -A_{\bf T})^{-1}v=\int_0^{+\infty}e^{-\lambda t}T(t)vdt \ \ \ \ (v \in \dt),\end{equation} 

which defines a quasimultiplier on $\iit$ if we apply it to some $v \in \Omega(\iit)\cap \dt.$ The fact that this quasimultiplier is regular is not completely obvious but follows from the previous discussion since $(\lambda I -A_{\bf T})^{-1}\in \jjt \subset {\mathcal QM}_r(\iit).$ Notice that since $\cup_{t>0}T(t)\iit$ is dense in $(\dt, \Vert . \Vert_{\dt}),$ $(\lambda I -A_{\bf T})^{-1}$ is characterized by the simpler formula

\begin{equation}(\lambda I -A_{\bf T})^{-1}T(s)v=e^{\lambda s}\int_s^{+\infty}e^{-\lambda t}T(t)vdt \ \ \ \ (s>0, v \in \iit).\end{equation} 

Using the correspondence between quasimultipliers on $\iit$ and partially defined operators on $X_1=\left [ \cup_{t>0}T(t)X\right ]^-$ and $X_0=X_1/\cap_{t>0}Ker(T(t)),$ we obtain a new approach to the resolvent formula and the other results of section 4.

\section{Holomorphic functional calculus}

For $\alpha \in \R,$ set $\Pi^+_{\alpha}:=\{ z \in \C \ : \ Re(z)>\alpha\},$ and set $\overline \Pi^+_{\alpha}:=\{ z \in \C \ : \ Re(z)\ge \alpha\}.$ Denote as usual by 
$H^1(\Pi^+_{\alpha})$ the space of all holomorphic functions $F$ on $\Pi^+_{\alpha}$ such that $\sup_{t>\alpha}\int_{-\infty}^{+\infty}\left |F(t+iy)\right |dy<+\infty,$ and set $\Vert F\Vert _1:=\sup_{t>0} \int_{-\infty}^{+\infty}\left |F(t+iy)\right |dy.$ Then   $(H^1(\Pi^+_{\alpha}), \Vert .\Vert_1)$ is a Banach space, the formula $F^*(\alpha +iy)=\lim_{t\to 0^+}F(\alpha +iy+t)$ defines a.e. a function $F^*$ on $\alpha +i\R,$ $\int_{-\infty}^{+\infty}\vert F^*(\alpha+iy)\vert dy<+\infty$, and we have

\begin{equation} \lim_{t\to 0^+}\int_{-\infty}^{+\infty} \vert F^*(\alpha +iy)-F(\alpha+t+iy)\vert dt=0.\end{equation} 

Moreover $\Vert F\Vert_1=\int_{-\infty}^{+\infty}\vert F^*(\alpha +iy)\vert dy,$ and if we set $f(t):={1\over 2\pi}e^{\alpha t}\int_{-\infty}^{+\infty}F^*(\alpha +iy)e^{iyt}dy,$
then it follows from the F. and M. Riesz theorem for the half-plane that $f(t)=0$ for $t<0,$ and ${\mathcal L}(f)(z)=F(z)$ for $z \in \Pi_{\alpha}^+.$ So we can define the inverse Laplace transform ${\mathcal L}^{-1}(F)$ of $F$ by the formula

\begin{equation} {\mathcal L}^{-1}(F)(t)={1\over 2\pi}e^{\alpha t}\int_{-\infty}^{+\infty}F^*(\alpha +iy)e^{iyt}dy.\end{equation}

It follows then from the inversion formula for Fourier transforms that ${\mathcal L}^{-1}(F)(t)={1\over 2\pi}e^{\beta t}\int_{-\infty}^{+\infty}F(\beta +iy)e^{iyt}dy$ for every $\beta >\alpha,$
and the tautological formula ${\mathcal L}^{-1}(F)={\mathcal L}^{-1}(F_{|_{\Pi_{\beta}^+}})$ holds for $F \in H^1(\Pi_{\alpha}^+), \beta > \alpha.$
All these results are well-known, see for example \cite{ri}, ch. 2.

Notice that if $F \in H^1(\Pi_{\alpha}^+),$ we have $\vert F(z)\vert =\left \vert \int_0^{+\infty}{\mathcal L}^{-1}(f)(t)e^{-zt}dt\right \vert \le {1\over 2\pi}\Vert F\Vert_1\int_0^{+\infty}e^{(\alpha -\beta)t}dt={1\over 2\pi}{\Vert F\Vert_1\over \beta -\alpha}$ for $Re(z)\ge \beta, \beta >\alpha.$ Hence $F_{|_{\Pi_{\beta}^+}}\in H^{\infty}(\Pi_{\beta}^+)$ for $\beta >\alpha,$ where $H^{\infty}(\Pi_{\beta}^+)$ denotes the algebra of all bounded holomorphic functions on $\Pi_\beta^+.$ This shows that $FG \in H^{1}(\Pi_{\gamma}^+)$ if $F\in H^1(\Pi_{\alpha}^+), G \in H^1(\Pi_{\beta}^+), \gamma >\beta \ge \alpha,$ and ${\mathcal L}({\mathcal L}^{-1}(F)*{\mathcal L}^{-1}(G))(z)= F(z)G(z)$ for $Re(z)\ge \gamma,$ which shows that ${\mathcal L}^{-1}(FG)={\mathcal L}^{-1}(F)*{\mathcal L}^{-1}(G).$ So $ \cup_{\alpha < \delta}H^1(\Pi_{\alpha}^+) \ \mbox{is an algebra for} \ \delta \in \R \cup{+\infty},$ and

\begin{equation} {\mathcal L}^{-1}(FG)={\mathcal L}^{-1}(F)*{\mathcal L}^{-1}(G) \ \ (F \in \cup_{\alpha < \delta}H^1(\Pi_{\alpha}^+), G \in \cup_{\alpha < \delta}H^1(\Pi_{\alpha}^+)).\end{equation}

\begin{prop}  Set, for $F \in H^1(\Pi_{\alpha}^+),$ $\alpha \in (-\infty, - log(\rho_{\bf T})),$

$$F(-A_{\bf T})=-{1\over 2\pi}\int_{-\infty}^{+\infty}F^*(\alpha +iy)(A_{\bf T} + (\alpha +iy)I)^{-1}dy \in \jjt \subset {\mathcal QM}_r(\iit).$$

Then $\mathcal L^{-1}(F)\in L^1_{\tilde \omega_{\bf T}},$ $F(-A_{\bf T})=\tilde \phi_{\bf T}((\mathcal L^{-1}(F))=F_{|_{\Pi_{\beta}}}(-A_{\bf T})=-{1\over 2\pi}\int_{-\infty}^{+\infty}F(\beta +iy)(A_{\bf T} + (\beta +iy)I)^{-1}dy$ for $\alpha < \beta <-log(\rho_{\bf T}),$
$\tilde \chi (F(-A_{\bf T}))=F(-\tilde \chi (A_{\bf T}))$ for $\chi \in \sia,$ and $FG(-A_{\bf T})= F(-A_{\bf T})G(-A_{\bf T})$ for  $F\in \cup_{\alpha < -log(\rho_{\bf T})}H^1(\Pi_{\alpha}^+),$
$G\in \cup_{\alpha < -log(\rho_{\bf T})}H^1(\Pi_{\alpha}^+).$
\end{prop}

Proof: We have $\vert {\mathcal L}^{-1}(F)(t)\vert \le e^{\alpha t}\Vert F\Vert_1$ for $t\ge 0,$ and so $\int_0^{+\infty}\vert {\mathcal L}^{-1}(F)(t)\vert \tilde \omega (t)dt\le \Vert F\Vert_1\int_0^{+\infty}e^{\alpha t}\tilde \omega_{\bf T}(t)dt<+\infty,$ since $\lim_{t\to +\infty}\tilde \omega(t)^{1/t}\le  \lim_{t\to +\infty} \omega(t)^{1/t}=\rho_{\bf T}.$ So ${\mathcal L}^{-1}(F)\in L^1_{\tilde \omega_{\bf T}}.$

Let $u \in \jjt.$ It follows from Fubini's theorem that we have

$$\tilde \phi_{\bf T}({\mathcal L}^{-1}(F))(u)={1\over 2\pi}\int_{0}^{+\infty}\left [\int_{-\infty}^{+\infty}e^{iyt}F^*(\alpha +iy)dy\right ]e^{\alpha t}T(t)udt$$
$$={1\over 2\pi}\int_{-\infty}^{+\infty}F^*(\alpha +iy)\left [ \int_0^{+\infty}e^{(\alpha +iy)t}T(t)udt\right ]dy$$ $$=-{1\over 2\pi}\int_{-\infty}^{+\infty}F^*(\alpha +iy)(A_{\bf T} +(\alpha +iy)I)^{-1}udy=F(-A_{\bf T})u.$$ Hence $F(-A_{\bf T})= \tilde \phi_{\bf T}((\mathcal L^{-1}(F))= \tilde \phi_{\bf T}({\mathcal L}^{-1}(F_{|_{\Pi_{\beta}}}))=F_{|_{\Pi_{\beta}}}(-A_{\bf T})=-{1\over 2\pi}\int_{-\infty}^{+\infty}F(\beta +iy)(A_{\bf T} + (\beta +iy)I)^{-1}dy$ for $\alpha < \beta <-log(\rho_{\bf T}).$

It follows from formula (15) that we have, for $\chi \in \sia,$

$$\tilde \chi (F(-A_{\bf T}))=\tilde \chi (\tilde \phi_{\bf T} ({\mathcal L}^{-1}(F))={\mathcal L}({\mathcal L}^{-1}(F))(-\tilde \chi(A_{\bf T}))=F(-\tilde \chi(A_{\bf T})).$$

Now let $F\in \cup_{\alpha < -log(\rho_{\bf T})}H^1(\Pi_{\alpha}^+),$
$G\in \cup_{\alpha < -log(\rho_{\bf T})}H^1(\Pi_{\alpha}^+).$ We have

$$FG(-A_{\bf T})=\tilde \phi_{\bf T} \left ({\mathcal L}^{-1}(FG)\right )=\tilde \phi_{\bf T} \left ( {\mathcal L}^{-1}(F)*{\mathcal L}^{-1}(G)\right)=\tilde \phi_{\bf T} \left ( {\mathcal L}^{-1}(F)\right )\tilde \phi _{\bf T} \left ({\mathcal L}^{-1}(G)\right )$$ $$=F(-A_{\bf T})G(-A_{\bf T}).$$ $\square$


Recall that a holomorphic function $F$ on $\Pi_{\alpha}^+$ is said to be outer if there exists a real-valued measurable function $u$ on $\R$ satisfying $\int_{-\infty}^{+\infty}{\vert u(t)\vert\over 1+t^2}dt<+\infty$ such that $F(z)= exp\left ({1\over \pi}\int_{-\infty}^{+\infty}{1+it(z-\alpha)\over it +z-\alpha}{u(t)\over 1+t^2}dt\right )$ for $z \in \Pi_{\alpha}^+,$ and in this situation we have $log\vert F^*(\alpha +iy)\vert =u(y)$ a.e., where $F^*(\alpha +iy)$ is defined for almost every $y\in \R$ by the formula $F^*(\alpha +iy)=\lim_{t\to 0^+}F(\alpha +t+iy),$   see for example \cite{roro}. The Smirnov class  ${\mathcal N}^+(\Pi_{\alpha}^+ )$ is the space of holomorphic functions $F$ on $\Pi_{\alpha}^+$ which can be written in the form $F=G/H$ where $G\in H^{\infty}(\Pi_{\alpha}^+)$ and $H \in H^{\infty}(\Pi_{\alpha}^+)$ is outer, so that the function $F^*: \alpha +iy \to \lim_{t\to 0^+}F(\alpha +iy +t)$ is defined a.e. on $\alpha +i\R.$
Since the restriction to $\Pi_{\beta}^+$ of an outer function on $\Pi_\alpha^+$ is an outer function on $\Pi_{\beta}^+$, we have $F_{|_{\Pi_{\beta}^+}}\in {\mathcal N}^+(\Pi_{\beta} ^+)$ for $F \in {\mathcal N}^+(\Pi_{\alpha}^+ )$ if $\beta >\alpha.$

The following property is a standard consequence of the transfer to the weighted convolution algebras $L^1_{u_\lambda}$ of Nyman's characterization of dense principal ideals of $L^1(\R^+),$ see \cite{ny}, but we prove it by a direct argument for the convenience of the reader.

\begin{lem} Let $\alpha < -log(\rho_{\bf T}),$ and let $F \in H^1(\Pi_\alpha^+)$ be outer. Then $F(-A_{\bf T})\in \Omega(\jjt).$
\end{lem}

Proof: The formula $u_n(t)=\min(-log(\vert F^*(\alpha +it)\vert), n)$ defines a.e. on $\R$ a function $u_n$ such that $\int_{-\infty}^{+\infty}{u_n(t)\over 1+t^2}dt <+\infty.$  Set, for $z \in \Pi_{\alpha}^+,$

$$F_n(z)= exp\left ({1\over \pi}\int_{-\infty}^{+\infty}{1+it(z-\alpha)\over it +z-\alpha}{u_n(t)\over 1+t^2}dt\right ).$$

Then $F_n(z)\le e^n,$ and $\vert F(z)F_n(z)\vert \le 1$ for $z \in \Pi_{\alpha}^+.$ It follows from the dominated convergence theorem  that we have

$$\lim_{n\to +\infty}\int_{-\infty}^{+\infty}\frac {\vert log \vert F^*_n(\alpha +it)\vert+u_n(t)\vert}{1+t^2}dt=\lim_{n\to +\infty}\int_{-\infty}^{+\infty}\frac { -log \vert F^*_n(\alpha +it)\vert-u_n(t)}{1+t^2}dt
=0,$$

and so $F(z)F_n(z)$ converges uniformly to 1 on $\Pi_{\beta}$ for every $\beta >\alpha.$ Let $\beta \in (\alpha, -log(\rho_{\bf T})),$ and let $G \in H^1(\Pi_{\alpha}).$ Then we have

$$\lim_{n\to +\infty}\int_{-\infty}^{+\infty} \vert G(\beta +it)- F(\beta +it)F_n(\beta+it)G(\beta+it)\vert dt=0,$$

and so $\lim_{n\to +\infty} \Vert {\mathcal L}^{-1}(G) - {\mathcal L}^{-1}(F)*{\mathcal L}^{-1}(F_nG)\Vert_{L^{\infty}(\R^+, e^{\beta t})}=0.$ Since $\tilde \phi_{\bf T}:L^{\infty}(\R^+, e^{\beta t})\to \jjt$ is continuous, this shows that $F(-A_{\bf T})\jjt$ contains $G(-A_{\bf T})$ for every $G \in H^1(\Pi_{\alpha}).$ Since ${\mathcal L}(f)\in \cap_{\alpha \in \R}H^1(\Pi_{\alpha})$ for every $f \in {\mathcal C}_0^{\infty}((0,+\infty)),$ this shows that $F(-A_{\bf t})\jjt$ is dense in $\jjt.$ $\square$

Let $F=G/H\in {\mathcal N}^+(\Pi_{\alpha} ),$ where $G\in H^{\infty}(\Pi_{\alpha})$ and $H \in H^{\infty}(\Pi_{\alpha})$ is outer, and let $U\in H^1(\Pi_{\alpha})\cap H^{\infty}(\Pi_{\alpha})$ be outer (for example one can take $U(z)= {1\over ( z-\alpha +1)^2}$). Then $F=GU/HU,$ where $GU \in H^1(\Pi_{\alpha})\cap H^{\infty}(\Pi_{\alpha})$ and where 
$HU \in H^1(\Pi_{\alpha})\cap H^{\infty}(\Pi_{\alpha})$ is outer. This observation allows to extend the holomorphic functional calculus to the family $\cup_{\alpha < -log(\rho_{\bf T})}
{\mathcal N}^+(\Pi_{\alpha} ).$

\begin{thm} For $F\in {\mathcal N}^+(\Pi_{\alpha}),$ $\alpha < -log(\rho_{\bf T}),$ set

$$F(-A_{\bf T})=(FH)(-A_{\bf T})/H(-A_{\bf T}) \in {\mathcal {QM}}(\jjt)={\mathcal {QM}}(\iit), $$ where $ H \in H^1(\Pi_{\alpha})$ is an outer function such that $FH \in H^1(\Pi_{\alpha}).$

(i)  The definition of $F(-A_{\bf T})$ does not depend on the choice of $H,$ $F(-A_{\bf T})=F_{|_{\Pi_\beta}}(-A_{\bf T})$ for $\beta \in (\alpha, -log(\rho_{\bf T})),$ $\tilde \chi (F(-A_{\bf T}))=F(-\tilde \chi (A_{\bf T}))$ for $\chi \in \sia,$ and $F(-A_{\bf T})\in {\mathcal {QM}_r}(\jjt)={\mathcal {QM}}_r({\iit})$ if, further, $F\in H^{\infty}(\Pi_{\alpha}).$

(ii)  $(FG)(-A_{\bf T})=F(-A_{\bf T})G(-A_{\bf T})$ for $F\in \cup_{\alpha < -log(\rho_{\bf T})}{\mathcal N}^+(\Pi_{\alpha}^+),$
$G\in \cup_{\alpha < -log(\rho_{\bf T})}{\mathcal N}^+(\Pi_{\alpha}^+).$ 

(iii) If $\mu$ is a measure on $[0,+\infty)$ such that $\int_0^{+\infty}e^{-\alpha t}d\vert \mu \vert (t)<+\infty$ for some $\alpha < -log(\rho_{\bf T}),$ then ${\mathcal L}(\mu)(-A_{\bf T})=\tilde \phi_{\bf T}(\mu).$ In particular $F(z)=e^{-tz}$ for some $t\ge 0,$ then $F(-A_{\bf T})=T(t).$

(iv) If $F(z)=-z,$ then $F(-A_{\bf T})=A_{\bf T}.$

\end{thm}
  
  Proof: It follows from the lemma that $H(-A_{\bf T})\in \Omega(\jjt)$ if $H \in H^1(\Pi_{\alpha})$ is outer, and so the formula $F(-A_{\bf T})= (FH)(-A_{\bf T})/H(-A_{\bf T})$ defines an element of ${\mathcal {QM}}(\jjt)={\mathcal {QM}}({\iit})$ if $ H \in H^1(\Pi_{\alpha})$ is an outer function such that $FH \in H^1(\Pi_{\alpha}).$ Let $H_1 \in H^1(\Pi_{\alpha})$ be another outer function such that $FH_1 \in H^1(\Pi_{\alpha}).$ It follows from the proposition that we have
  
  $$(FH)(-A_{\bf T})H_1(-A_{\bf T})=(FHH_1)(-A_{\bf T})=(FH_1)(-A_{\bf T})H(-A_{\bf T}),$$
  
  and so the definition of $F(-A_{\bf T})$ does not depend on the choice of the outer function $H \in H^1(\Pi_{\alpha}).$
  
  Let $\beta \in (\alpha, -log(\rho_{\bf T})).$ Then $(FH)_{|_{\Pi_{\beta}}}\in H^1(\Pi_{\beta}), H_{|_{\Pi_{\beta}}}\in H^1(\Pi_{\beta})$ is outer,  and it follows from the proposition that we have
  
  $$F_{|_{\Pi_{\beta}}}(-A_{\bf T})=(FH)_{|_{\Pi_{\beta}}}(-A_{\bf T})/H_{|_{\Pi_{\beta}}}(-A_{\bf T})=(FH)(-A_{\bf T})/H(-A_{\bf T})=F(-A_{\bf T}).$$
  
  Let $\chi \in \sia.$ It follows from the definition of $\tilde \chi$ and from the proposition that we have
  
  $$\tilde \chi (F(-A_{\bf T}))$$ $$=\tilde \chi ((FH)(-A_{\bf T}))/\tilde \chi (H(-A_{\bf T}))=(FH)(-\tilde \chi(A_{\bf T}))/H(-\tilde \chi(A_{\bf T}))= F(-\tilde \chi (A_{\bf T})).$$

  Let $F_1\in {\mathcal N}^+(\Pi_{\alpha_1}^+),$ let $F_2\in {\mathcal N}^+(\Pi_{\alpha_2}^+),$ where $\alpha_1<-log(\rho_{\bf T}), \alpha_2<log(\rho_{\bf T}),$ let $H_1\in H^1(\Pi_{\alpha_1}^+)$ such that $F_1H_1\in H^1(\Pi_{\alpha_1}^+),$ let $H_2\in H^1(\Pi_{\alpha_2}^+)$ such that $F_2H_2\in H^1(\Pi_{\alpha_2}^+),$ and let $\beta \in (\max(\alpha_1,\alpha_2), -log(\rho_{\bf T}).$ Denote respectively by $\tilde F_1, \tilde F_2, \tilde H_1$ and $\tilde H_2$ the restrictions of   $ F_1,  F_2,  H_1$ and $ H_2$ to $\Pi_{\beta}^+.$ Then $\tilde H_1\tilde H_2 \in H^1(\Pi_{\beta}^+)\cap H^{\infty}({\Pi_\beta}^+)$ is outer, $\tilde H_1\tilde H_2\tilde F_1\tilde F_2 \in H^1(\Pi_{\beta}^+)\cap H^{\infty}(\Pi_{\beta}^+),$ and we have
  
  $$(F_1F_2)(-A_{\bf T})=(\tilde{F_1}\tilde {F_2})(-A_{\bf T})=(\tilde F_1 \tilde F_2 \tilde H_1\tilde H_2)(-A_{\bf T})/(\tilde H_1\tilde H_2)(-A_{\bf T})$$ 
  
$$= \left ( (\tilde F_1\tilde H_1)(-A_{\bf T})(\tilde F_2\tilde H_2)(-A_{\bf T})\right )/\left (\tilde H_1(-A_{\bf T})\tilde H_1(-A_{\bf T})\right )=\tilde F_1(-A_{\bf T})\tilde F_2(-A_{\bf T})$$
$$= F_1(-A_{\bf T}) F_2(-A_{\bf T}),$$ which proves (ii).
     
  Let  $F\in H^{\infty}(\Pi_{\alpha}^+)$ for some $\alpha < -log(\rho_{\bf T}),$ and let $H\in H^1(\Pi_{\alpha}^+)$ be outer. We have, for $n\ge 1,$
  
 $$\Vert H(-A_{\bf T})(F(-A_{\bf T})^n\Vert_{\jjt}= \Vert (HF^n)(-A_{\bf T})\Vert_{\jjt}$$ $$\le \left [\sup_{t\ge 0}e^{-\alpha t}{\mathcal L}^{-1}(HF^n)(t)\right ]\int_0^{+\infty}e^{\alpha t} \Vert T(t)\Vert dt$$ $$\le \left [{1\over 2\pi}\int_{-\infty}^{+\infty}\vert H(\alpha +iy)\vert dy\int _0^{+\infty}e^{\alpha t} \Vert T(t)\Vert dt\right ]\Vert F\Vert^n_{H^{\infty}(\Pi_{\alpha}^+)}.$$
 
 This shows that $F(-A_{\bf T})$ is a regular mutiplier on $\jjt,$ and $F(-A_{\bf T})\in {\mathcal {QM}}_r(\jjt)={\mathcal {QM}}_r(\iit),$ which completes the proof of (i).
 
Now let $\alpha < -log(\rho_{\bf T}),$ let $\mu$ be a measure on $[0,+\infty)$ such that $\int_0^{+\infty}e^{\alpha s}d\vert \mu\vert (s)<+\infty$ for some $\alpha < -log(\rho_{\bf T}),$ and let $F \in H^1(\Pi_{\alpha}^+)$ be an outer function. We have

$${\mathcal L}(\mu)(-A_{\bf T})=({\mathcal L}(\mu)F)(-A_{\bf T})/F(-A_{\bf T})= \tilde \phi_{\bf T} (\mu *{\mathcal L}^{-1}(F))/\tilde \phi _{\bf T}({\mathcal L}^{-1}(F))=\tilde \phi_{\bf T} (\mu).$$

In particular if $F(z)=e^{-tz}$ for some $t \ge 0$ we have $F(-A_{\bf T})={\mathcal L}(\delta_t)(-A_{\bf T})=\tilde \phi_{\bf T}(\delta_t)=T(t),$ which completes the proof of $(iii).$

Let again $\alpha < -log(\rho_{\bf T}),$ and set $f(t)=t^2e^{(\alpha -1)t}$ for $t\ge 0.$ Then ${\mathcal L}(f)(z)={2\over (z+1-\alpha)^3}$ for $z \in {\overline \Pi} _{\alpha}^+,$ so that ${\mathcal L}\in H^1(\Pi_{\alpha})$ is outer. Since $f \in {\mathcal C}^1([0,+\infty))\cap L^1_{\tilde \omega_{\bf T}},$ and since $f(0)=0,$ it follows from definition 8.3 that $A_{\bf T}= -\tilde \phi _{\bf T}(f')/\tilde \phi _{\bf T}(f).$ Set $F(z)=-z.$ We have $-{\mathcal L}(f')=F{\mathcal L}(f)\in H^1(\Pi_{\alpha}),$ and

$$F(-A_{\bf T})=-{\mathcal L}(f')(-A_{\bf T})/{\mathcal L}(f)(-A_{\bf T})=-\tilde \phi_{\bf T}(f')/\tilde \phi_{\bf T}(f)=A_{\bf T},$$

which proves the tautological assertion (iv). $\square$

Let $F \in \cup_{\alpha < -log(\rho_{\bf T})}{\mathcal N}^+(\Pi^+_{\alpha})$.The quasimultiplier $F(-A_{\bf T})=(FH)(-A_{\bf T})/H(-A_{\bf T}),$ where $ H \in H^1(\Pi_{\alpha})$ is an outer function such that $FH \in H^1(\Pi_{\alpha})$ defined above  is an element of ${\mathcal{QM}}(\jjt).$ In order to interpret $F(-A_{\bf T})$ as an element of  ${\mathcal{QM}}(\iit)$ we just have to write $F(-A_{\bf T})=u(FH)(-A_{\bf T})/uH(-A_{\bf T}),$ where $u \in \Omega(\dt).$ 

Set again $u_{\alpha}(t)=e^{\alpha t}$ and $\omega_{\alpha}(t)=e^{\alpha_t}sup_{s \ge t}e^{-\alpha s}\Vert T(s)\Vert.$ We can also write directly $F(-A_{\bf T})=(FH_1H_2)(-A_{\bf T})/H_1H_2(-A_{\bf T}),$ where $ H_1 \in H^1(\Pi_{\alpha})$ is an outer function such that $FH_1 \in H^1(\Pi_{\alpha})$ and where $H_2={\mathcal L}(f)$ with $f \in \Omega(L^1_{u_{\alpha}})\cap 
L^1(u_{\omega_{\alpha}}).$ In this case ${\mathcal L}(f)\in H^{\infty}(\Pi^+_{\alpha})$ is outer on $\Pi^+_{\alpha},$ and so $H_1H_2\in H^1(\Pi_{\alpha}^+)$ is outer, $FH_1H_2 \in H^1(\Pi_{\alpha}^+),$ $H_1(-A_{\bf T})\in \Omega(\jjt), H_2(-A_{\bf T})=\tilde \phi_{\bf T}(f ) \in \tilde \theta_{\bf T}\left (\Omega(L^1_{\omega_ \alpha})\right )\subset \Omega(\dt),$ so that $H_1H_2(-A_{\bf T})=H_1(-A_{\bf T})H_2(-A_{\bf T})\in \Omega(\dt)\subset \Omega(\iit),$ and $(FH_1H_2)(-A_{\bf T})=(FH_1)(-A_{\bf T})H_2(-A_{\bf T})\in \dt \subset \iit,$ so that the formula $F(-A_{\bf T})=(FH_1H_2)(-A_{\bf T})/H_1H_2(-A_{\bf T})$ defines explicitely a quasimultiplier on $\iit.$

Notice that this functional calculus applies to every function $F \in \cup_{\alpha < -log(\rho_{\bf T})}{\mathcal L}(L^1_{u_{-\alpha}}),$ since ${\mathcal L}(f)\in H^{\infty}(\Pi_{\alpha}^+)\cap{\mathcal C}(\overline {\Pi}_{\alpha}^+)$ for $f \in L^1_{u_{-\alpha}},$ but in general $F(-A_{\bf T})$ will be written as a quotient. Set $v_{n,\alpha}(t)=(n-\alpha)^2te^{-(n -\alpha)t}$ for $t\ge 0.$ Then $\int_0^{+\infty}v_{n,\alpha}(t)dt=1,$ and since $\lim_{n\to +\infty}\int_{\delta}^{+\infty}e^{-\alpha t}v_{n,\alpha}(t)dt=0$ for every $\delta >0,$ $(v_{n,\alpha})_{n\ge 1}$ is a sequential bounded approximate identity for the convolution algebra $L^1_{\omega_{\bf T}},$ and $\lim_{n\to +\infty}\Vert \tilde  \phi_{\bf T}(f*v_{n, \alpha})-\phi_{\bf T}(f)\Vert _{\jjt}=0.$ Since ${\mathcal L}(v_{n,\alpha})(z)={(n-\alpha))^2\over (z+n-\alpha)^2}$ for $Re(z) \ge \alpha,$ we obtain, for $f \in L^1_{u_{-\alpha}}, \alpha <-log(\rho_{\bf T}),$

\begin{equation}\lim_{n\to +\infty} \left \Vert {\mathcal L}(f)(-A_{\bf T}) +{1\over 2\pi}\int_{-\infty}^{+\infty}{(n-\alpha)^2\over (z+n-\alpha)^2}{\mathcal L}(F)(\alpha +iy)(A_{\bf T} +(\alpha +iy)I)^{-1}dy\right \Vert_{\jjt}=0. \end{equation}

When the semigroup ${\bf T}=(T(t))_{t>0}$ is quasinilpotent the class $\cup_{\alpha \in  \R}{\mathcal N}^+(\Pi_{\alpha}^+)$ seems to be the largest class on which the functional calculus $F \to F(-A_{\bf T})$ can be defined without other hypothesis on the Arveson weakly continous semigroup ${\bf T}.$ When the semigroup is not quasinilpotent, the convolution algebra $L^1_{\tilde \omega_{\bf T}}$ contains strictly the algebra $\cup_{\alpha <-log(\rho_{\bf T})}L^1_{u_{-\alpha}},$ but we can still write, for  $f\in L^1_{\tilde \omega_{\bf T}},$ the limit being taken with respect to the norm $\Vert . \Vert _{\jjt},$

$$ {\mathcal L}(f)(-A_{\bf T})=\tilde \phi_{\bf T}(f)=\lim_{\epsilon \to 0^+}\tilde \phi_{\bf T}(fu_{-\epsilon})={\mathcal L}(fu_{-\epsilon})(-A_{\bf T}),$$

and ${\mathcal L}(fu_{-\epsilon})(z) = {\mathcal L}(f)(z+\epsilon).$ The quasimultiplier ${\mathcal L}(fu_{-\epsilon})(-A_{\bf T})$ may be then computed by using (21), with $\alpha =-\rho_{\bf T} -\epsilon.$

When the semigroup is analytic on a sector, one can use Cauchy's formula to define $F(-A_{\bf T})$ for functions $F$ belonging to the $H^1$ space of a suitable family of sectors. We will not do it here.

IMB, UMR 5251

Universit\'e de Bordeaux

351, cours de la Lib\'eration

33405 - Talence

esterle@math.u-bordeaux1.fr
  

\begin{thebibliography}{55}
\bibitem{ri} F. Ricci, Hardy spaces in one variable, 
 
 {\small{\color{blue}homepage.sns.it/fricci/papers/hardy.pdf}}
\end{thebibliography}
  \end{document}